\documentclass[11pt]{article}
\usepackage{epsfig,epic,eepic,amsmath}
\input amssym.def
\input amssym.tex

\newcommand{\qed}{\rule{3mm}{3mm}}
\newcommand{\itbf}{\itshape\bfseries}

\newcommand{\aaa}{{\frak a}}

\newcommand{\ee}{{\frak e}}

\newcommand{\be}{{\bf e}}
\newcommand{\bn}{{\bf n}}
\newcommand{\ba}{{\bf a}}
\newcommand{\bb}{{\bf b}}
\newcommand{\bO}{{\boldsymbol 0}}

\newcommand{\beps}{{\boldsymbol\epsilon}}

\newcommand{\cD}{{\cal D}}
\newcommand{\cG}{{\cal G}}
\newcommand{\cH}{{\cal H}}
\newcommand{\cI}{{\cal I}}

\newcommand{\cS}{{\cal S}}

\newcommand{\bD}{{\bf D}}
\newcommand{\bG}{{\bf G}}

\newcommand{\bbR}{{\Bbb R}}
\newcommand{\bbC}{{\Bbb C}}
\newcommand{\bbZ}{{\Bbb Z}}
\newcommand{\bbS}{{\Bbb S}}

\newtheorem{theorem}{Theorem}
\newtheorem{proposition}[theorem]{Proposition}
\newtheorem{lemma}[theorem]{Lemma}
\newtheorem{corollary}[theorem]{Corollary}
\newtheorem{definition}[theorem]{Definition}

\begin{document}

\title{Linear and nonlinear theories of discrete
analytic functions. \\ Integrable structure and isomonodromic
Green's function}

\author{
Alexander I. Bobenko
\thanks{Institut f\"ur Mathematik,
Fachbereich II, TU Berlin, Str. des 17. Juni 136, 10623 Berlin,
Germany. E--mail: {\tt bobenko@math.tu-berlin.de}. Partially supported 
by the DFG Research Center ``Mathematics for key technologies'' (FZT86) 
in Berlin.} \and
Christian Mercat
\thanks{D\'epartement de  Math\'ematiques, Universit\'e
Montpellier II, Place Eug\`ene Bataillon, 34095 Montpellier Cedex
5, France. E--mail: {\tt mercat@math.univ-montp2.fr}. Partially supported
by DFG in the frame of SFB288 ``Differential geometry and quantum physics''.} 
\and
Yuri B. Suris
\thanks{Institut f\"ur Mathematik, Fachbereich II, TU Berlin,
Str. des 17. Juni 136, 10623 Berlin, Germany. E--mail: {\tt
suris@sfb288.math.tu-berlin.de}. Supported 
by the DFG Research Center ``Mathematics for key technologies'' (FZT86) 
in Berlin.}}

\maketitle

{\small {\bf Abstract.} Two discretizations, linear and nonlinear, of 
basic notions of the complex analysis are considered. The underlying lattice 
is an arbitrary quasicrystallic rhombic tiling of a plane. The linear theory 
is based on the
discrete Cauchy-Riemann equations, the nonlinear one is based on the notion
of circle patterns. We clarify the role of the rhombic condition in both
theories: under this condition the corresponding equations are integrable
(in the sense of 3D consistency, which yields also the existense of zero
curvature representations, B\"acklund transformations etc.). We demonstrate
that in some precise sense the linear theory is a linearization of the 
nonlinear one: the tangent space to a set of integrable circle patterns at an
isoradial point consists of discrete holomorphic functions which
take real (imaginary) values on two sublattices. 
We extend solutions of the basic equations of both theories
to $\bbZ^d$, where $d$ is the number of different edge slopes of the
quasicrystallic tiling.
In the linear theory, we give an integral representation of an arbitrary
discrete holomorphic function, thus proving the density of discrete 
exponential functions. We introduce the $d$-dimensional discrete logarithmic 
function which is a generalization of Kenyon's discrete Green's function, and 
uncover several new properties of this function. We prove that it is 
an isomonodromic solution of the discrete Cauchy-Riemann equations, and
that it is a tangent vector to the space of integrable
circle patterns along the family of isomonodromic discrete power functions. 
 }

\section{Introduction}
\label{Sect intro}

There is currently much interest in finding discrete counterparts
of various structures of the classical (continuous, smooth)
mathematics. In the present paper we are dealing with the
discretization of the classical complex analysis.

There are two approaches to this problem. The first one, which we
shall call the {\it linear theory}, is based on a discretization
of the Cauchy-Riemann equations. Since the latter are linear,
straightforward discretizations are linear as well. A
discretization preserving apparently the most number of important
structural features has been developed in \cite{F, D1, D2, M1, K}.
The first two references are dealing with {\it discrete
holomorphic functions} $f:\bbZ^2\to\bbC$ on the regular square
lattice, satisfying the following {\it discrete Cauchy-Riemann
equations}:
\begin{equation}\label{CR intro}
f_{m,n+1}-f_{m+1,n}=i(f_{m+1,n+1}-f_{m,n}).
\end{equation}
A pioneering step was undertaken by Duffin \cite{D2}, where the
combinatorics of $\bbZ^2$ was given up in favor of arbitrary
planar graphs with rhombic faces. A far reaching generalization of
these ideas is given in \cite{M1}, where the linear theory is extended 
to discrete Riemann surfaces. Planar graphs with rhombic
faces are called {\it critical} in \cite{M1}. Kenyon \cite{K}
developed a theory of the Dirac operator and constructed Green's
function in the framework of the linear theory on critical graphs.
See \cite{CY, G} for combinatorial, resp. numerical aspects of Green's
functions on graphs.

The second approach, which we will call the {\it nonlinear
theory}, is based on the ideas by Thurston \cite{T}, and
declares {\it circle patterns} to be natural discrete analogs of
analytic functions \cite{BeS, DS, Sch, S}. One of the most important
achievements of this theory is the proof that the
Riemann map can be (constructively) approximated by circle
packings \cite{RS, MR, HS}. The variational approach to circle
patterns is discussed in detail in \cite{BSp}. The word
``nonlinear'' refers to the basic feature of equations describing
circle patterns. Often, the so-called {\em cross-ratio system} is
used for this. For a function $f:\bbZ^2\to\bbC$ on the regular
square lattice, this system was introduced in \cite{NC}:
\begin{equation}\label{cr intro}
\frac{(f_{m+1,n}-f_{m,n})(f_{m+1,n+1}-f_{m,n+1})}
{(f_{m,n+1}-f_{m,n})(f_{m+1,n+1}-f_{m+1,n})}=-1.
\end{equation}
For circle patterns with more sophisticated combinatorics, a
generalization of this system to an arbitrary quad-graph (planar
graph with quadrilateral faces) is required \cite{BS}.

It is not difficult to see in what sense solutions of equations
like (\ref{CR intro}), (\ref{cr intro}) can be considered as
discretized analytic functions. Indeed, assume that $\bbZ^2$ is
embedded in the complex plane $\bbC$ with the grid size
$\varepsilon$, i.e., the pair $(m,n)\in\bbZ^2$ corresponds to
$(m+in)\varepsilon\in\bbC$. Then restrictions of analytic
functions to this grid satisfy the corresponding equations up to
$O(\varepsilon^2)$. More precisely, if $f:\bbC\to\bbC$ is
analytic, then
\[
\frac{f(z+i\varepsilon)-f(z+\varepsilon)}
{f(z+\varepsilon+i\varepsilon)-f(z)}=i+O(\varepsilon^2),
\]
and
\[
\frac{\big(f(z+\varepsilon)-f(z)\big)
      \big(f(z+\varepsilon+i\varepsilon)-f(z+i\varepsilon)\big)}
     {\big(f(z+i\varepsilon)-f(z)\big)
      \big(f(z+\varepsilon+i\varepsilon)-f(z+\varepsilon)\big)}=
   -1+O(\varepsilon^2).
\]
Similar relations hold on more general graphs.

For a long time, the linear and the nonlinear theories of discrete
analytic functions were considered separately. In the present
paper, we show that in some precise sense the former is a {\it
linearization} of the latter. We work in the set-up of rhombic 
tilings of a plane. The theory becomes especially rich for 
{\em quasicrystallic} tilings, -- those with a finite 
number of different edge slopes. This class includes double periodic
tilings (which are naturally considered on a torus), as well as
non-periodic ones, like the Penrose tiling. We clarify the
importance of rhombic embeddings of quad-graphs in both the linear and 
the nonlinear theories. Namely, we show that the rhombic property implies 
(actually, is almost synonymous with) {\it integrability}. Note that
interrelations of circle patterns with the theory of integrable
systems were already uncovered and studied in \cite{BP, AB1, AB2,
BHS, BH}. Note also that some of the ideas behind our unified
treatment of integrability of linear and nonlinear systems, such
as the use of zero curvature representations in both situations,
are similar to the philosophy of Fokas's unified transform method for
linear and nonlinear differential equations based on the Riemann-Hilbert
boundary problem \cite{Fo}. Our main results are the following.
\begin{itemize}
\item Discrete Cauchy-Riemann equations on a rhombically embedded quad-graph
 $\cD$, with weights given by quotients of diagonals of the corresponding
 rhombi, are integrable. Integrability is understood here as 3D consistency
 \cite{BS}. Therefore, discrete holomorphic functions on rhombic
 embeddings can (and should) be extended to multidimensional lattices.
 In particular, discrete holomorphic functions on a quasicrystallic rhombic
 embedding $\cD$ with $d$ different edge slopes can be considered as 
 restrictions of discrete holomorphic functions on $\bbZ^d$ to a 
 certain two-dimensional subcomplex $\Omega_\cD$ in $\bbZ^d$.
\item  Cross-ratio equations on a rhombically embedded quad-graph $\cD$,
 with cross-ratios read off the corresponding rhombi, are integrable as well.
 Therefore, solutions of the cross-ratio equations on a quasicrystallic
 rhombic embedding $\cD$ are naturally extended to $\bbZ^d$.
\item For a circle pattern, the centers and the intersection points of the
 circles yield a solution of cross-ratio equations, with the cross-ratios
 depending on the pairwise intersection angles of the circles. We say that
 a circle pattern is integrable, if the corresponding cross-ratio system is
 integrable. The combinatorics and intersection angles belong to an integrable
 circle pattern, if and only if they admit an isoradial realization.
 This latter realization gives a rhombic immersion of the corresponding
 quad-graph, and generates also a dual isoradial circle pattern.
 An integrable circle pattern can be alternatively described by the radii
 of the circles and the rotation angles of the configurations at the
 intersection points with respect to the isoradial realization. These data
 comprise a solution of an integrable {\em Hirota system}.
\item The tangent space to the set of integrable circle patterns, at the
 point corresponding to an isoradial pattern, coincides with the space
 of discrete holomorphic functions on the corresponding rhombically embedded
 quad-graph, which take real (resp. pure imaginary) values on the white
 (resp. black) vertices. This holds in the description of circle patterns in
 terms of circle radii and rotation angles at the intersection points (Hirota
 equations). Discrete holomorphic functions obtained from these ones by
 discrete integration, comprise the tangent space to the set of integrable 
 circle patterns, described in terms of circle centers and 
 intersection points (cross-ratio equations).
\item We define (in the linear theory) discrete exponential functions on
 $\bbZ^d$, and prove that they are dense in the space of discrete
 holomorphic functions, growing not faster than exponentially.
\item We define (in the linear theory) a discrete logarithmic function on
 $\bbZ^d$, or, better, on a branched covering of certain $d$-dimensional
 octants
    \footnote{We use this term for a subset of $\bbZ^d$ defined by
    fixing one of $2^d$ possible combinations of signs of the coordinates.
    An octant in the proper sense corresponds to $d=3$, while by $d=2$ this
    object is called quadrant.}
 $S_m\subset\bbZ^d$, $m=1,\ldots,2d$. On each such octant, the
 discrete logarithmic function is discrete holomorphic, with the
 distinctive property of being isomonodromic, in the sense of the integrable
 systems theory. We show that the real part of the discrete logarithmic
 function restricted to a surface $\Omega_\cD$ in $\bbZ^d$ coming from a
 quasicrystallic quad-graph $\cD$ is nothing but Green's function
 found in \cite{K}. The integral representation of Green's function given
 in \cite{K} is derived within the isomonodromic approach.
\item We define (in the nonlinear theory) discrete power functions
 $w^{\gamma-1}$ (resp. $z^\gamma$) on the same branched covering of octants
 $S_m\subset\bbZ^d$, $m=1,\ldots,2d$, where the discrete logarithmic function
 is defined. On each such sector, discrete $w^{\gamma-1}$ (resp. $z^\gamma$)
 is an isomonodromic solution of the Hirota (resp. cross-ratio) system.
 The tangent vector to the space of integrable circle patterns along 
 the curve consisting of patterns $w^{\gamma-1}$, at the isoradial point
 corresponding to $\gamma=1$, is shown to be the discrete logarithmic
 function.
\end{itemize}
In conclusion, we point out some generalizations of the concepts
and results of this paper for the non-rhombic case. 
\medskip

{\bf Acknowledgements.} Numerous discussions and collaboration with Boris 
Springborn were very important for this research. We thank also Tim Hoffmann,
Ulrich Pinkall and G\"unter Ziegler for discussions.

\section{\!\!Discrete harmonic and discrete holomorphic functions on graphs}
\label{Sect Laplace}

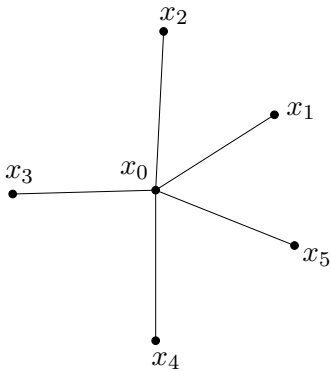
\begin{figure}[htbp]
    \setlength{\unitlength}{15pt}
\begin{center}
\begin{picture}(8,9)
\put(4,4){\circle*{0.2}} \put(4.2,8){\circle*{0.2}}
\put(7,5.9){\circle*{0.2}} \put(7.5,2.6){\circle*{0.2}}
\put(4,0.2){\circle*{0.2}} \put(0.4,3.9){\circle*{0.2}}
\path(4,4)(4.2,8) \path(4,4)(7,5.9) \path(4,4)(7.5,2.6)
\path(4,4)(4,0.2) \path(4,4)(0.4,3.9) \put(3.1,4.4){$x_0$}
\put(7.3,5.9){$x_1$} \put(4.1,8.3){$x_2$} \put(0.2,4.3){$x_3$}
\put(3.9,-0.4){$x_4$} \put(7.7,2.2){$x_5$}
\end{picture}
\end{center}
    \caption{The star of the vertex $x_0$ in the graph $\cG$ .}
        \label{star}
    \end{figure}

We denote by $V(\cG)$, $E(\cG)$ and $\vec{E}(\cG)$ the sets of vertices, 
undirected and directed edges of a graph $\cG$, respectively. 
Let there be given a complex-valued function
$\nu:E(\cG)\to\bbC$ on the edges. Then the {\it Laplacian}
$\Delta$ corresponding to the weight function $\nu$ is the
operator acting on functions $f:V(\cG)\to\bbC$ by
\begin{equation}\label{Laplace}
(\Delta f)(x_0)=\sum_{x\sim x_0}\nu(x_0,x)(f(x)-f(x_0)).
\end{equation}
Here the summation is extended over the set of all vertices $x$
connected to $x_0$ by an edge. We will use the notation ${\rm
star}(x_0)={\rm star}(x_0;\cG)$ for the set of all edges of $\cG$
incident to $x_0$, see Fig.\,\ref{star}.
\begin{definition}
A function $f:V(\cG)\to\bbC$ is called {\itbf discrete harmonic}
(with respect to the weights $\nu$), if $\Delta f=0$.
\end{definition}
Of course, the most interesting case of these notions is that of
real positive weights $\nu:E(\cG)\to\bbR_+$.
%------------------------------------------------------------------
%------------------------------------------------------------------------

In the continuous case, there is a canonical correspondence
between harmonic and holomorphic functions on $\bbC$: the real and
the imaginary parts of a holomorphic function are harmonic, and
any real-valued harmonic function can be considered as a real part
of a holomorphic function. This relation can be generalized for
functions on graphs, but these two classes of functions live then
on different graphs. Discrete holomorphic functions live on
quad-graphs.
\begin{definition}
A cell decomposition $\cD$ of the plane $\bbC$ is called a {\itbf
quad-graph}, if all its faces are quadrilaterals.
\end{definition}
A more general version of this definition deals with cell
decompositions of an arbitrary oriented surface. So, quad-graphs
are not just graphs, but are additionally assumed to be embedded
in an oriented surface; we will deal with the case of $\bbC$ only.

To establish a relation with discrete harmonic functions, we
consider the latter ones on graphs $\cG$ with an
additional structure, namely on those that come from general
(not necessarily quadrilateral) cell decompositions of $\bbC$. 
We will denote by $F(\cG)$ the set of faces (2-cells) of $\cG$. 
To any such $\cG$ there corresponds canonically a combinatorial
quad-graph called its double (or diamond \cite{M1}),
constructed from $\cG$ and its dual $\cG^*$. Recall that, in general, 
a {\it dual cell decomposition} $\cG^*$ is only defined up to isotopy, 
but it can be fixed uniquely with the help of the Voronoi/Delaunay 
construction. The dual $\cG^*$ is characterized as follows.
Vertices of $\cG^*$ are in a one-to-one correspondence to faces of
$\cG$, see Fig. \ref{dual vertex}. Each $\ee\in E(\cG)$ separates
two faces of $\cG$, which in turn correspond to two vertices of
$\cG^*$. It is declared that these two vertices are connected by
the edge $\ee^*\in E(\cG^*)$ dual to $\ee$. Finally,
the faces of $\cG^*$ are in a one-to-one correspondence with the
vertices of $\cG$: if $x_0\in V(\cG)$, and $x_1,\ldots,x_n\in
V(\cG)$ are its neighbors connected with $x_0$ by the edges
$\ee_1=(x_0,x_1),\ldots, \ee_n=(x_0,x_n)\in E(\cG)$, then the face
of $\cG^*$ corresponding to $x_0$ is defined by its boundary
$\ee_1^*\cup\ldots\cup\ee_n^*$ (cf. Fig. \ref{dual face}). 
If one assigns a direction to
an edge $\ee\in E(\cG)$, then it will be assumed that the dual
edge $\ee^*\in E(\cG^*)$ is also directed, in a way consistent
with the orientation of the underlying surface, namely so that the
pair $(\ee,\ee^*)$ is oriented directly at its crossing point.
This orientation convention implies that $\ee^{**}=-\ee$. 
%--------------------------------------------------------------------------
\begin{figure}[htbp]
    \setlength{\unitlength}{30pt}
    \begin{minipage}[t]{200pt}
\begin{picture}(6,5)
\put(3,1){\circle*{0.15}} \put(4.5,1.5){\circle*{0.15}}
\put(5,2.7){\circle*{0.15}} \put(3.5,4){\circle*{0.15}}
\put(2.1,2.5){\circle*{0.15}} \put(3.7,2.6){\circle{0.15}}
\path(3,1)(4.5,1.5)  \path(3,1)(2.8,0.6) \path(3,1)(3.2,0.6)
\path(4.5,1.5)(5,2.7) \path(4.5,1.5)(4.9,1.3) \path(5,2.7)(3.5,4)
\path(5,2.7)(5.3,2.4)  \path(5,2.7)(5.3,3) \path(3.5,4)(2.1,2.5)
\path(3.5,4)(3.5,4.4) \path(2.1,2.5)(3,1)  \path(2.1,2.5)(1.7,2.5)
\path(2.1,2.5)(1.8,2.3)  \path(2.1,2.5)(1.8,2.7)
\dashline[+30]{0.2}(3.7,2.6)(4.9,4.2)
\dashline[+30]{0.2}(3.7,2.6)(5.5,1.7)
\dashline[+30]{0.2}(3.7,2.6)(3.9,0.6)
\dashline[+30]{0.2}(3.7,2.6)(2.2,1.5)
\dashline[+30]{0.2}(3.7,2.6)(2.4,3.6)
\end{picture}
    \caption{Vertex of $\cG^*$ dual to a face of $\cG$.}
        \label{dual vertex}
    \end{minipage}\hfill
\begin{minipage}[t]{200pt}
    \setlength{\unitlength}{15pt}
\begin{picture}(7,9)(-2,0)
\put(4,4){\circle*{0.3}} \put(4.2,8){\circle*{0.3}}
\put(7,5.9){\circle*{0.3}} \put(7.5,2.6){\circle*{0.3}}
\put(4,0.2){\circle*{0.3}} \put(0.4,3.9){\circle*{0.3}}
\put(5,6){\circle{0.3}} \put(3,5.7){\circle{0.3}}
\put(3,2.1){\circle{0.3}} \put(5.1,2.5){\circle{0.3}}
\put(6,4.2){\circle{0.3}} \dashline[+30]{0.2}(5,6)(3,5.7)
\dashline[+30]{0.2}(3,5.7)(3,2.1)
\dashline[+30]{0.2}(3,2.1)(5.1,2.5)
\dashline[+30]{0.2}(5.1,2.5)(6,4.2)
\dashline[+30]{0.2}(6,4.2)(5,6) \path(4,4)(4.2,8)
\path(4,4)(7,5.9) \path(4,4)(7.5,2.6) \path(4,4)(4,0.2)
\path(4,4)(0.4,3.9) \put(4.2,3.3){$x_0$}\put(7.3,5.9){$x_1$}
\put(4.1,8.3){$x_2$}\put(0,4.3){$x_3$}
\put(3.9,-0.3){$x_4$}\put(7.5,2.1){$x_5$}
\end{picture}
    \caption{Face of $\cG^*$ dual to a vertex of $\cG$.}
        \label{dual face}
    \end{minipage}
\end{figure}
%----------------------------------------------------------------------

Now the {\it double} $\cD$ is constructed from $\cG$, $\cG^*$ as
follows. The set of vertices of the double $\cD$ is
$V(\cD)=V(\cG)\sqcup V(\cG^*)$. Each pair of dual edges, say 
$\ee=(x_0,x_1)\in E(\cG)$ and $\ee^*=(y_0,y_1)\in E(\cG^*)$, defines a
quadrilateral $(x_0,y_0,x_1,y_1)$. These quadrilaterals constitute the 
faces of the cell decomposition (quad-graph) $\cD$, see 
Fig.\,\ref{diamond}. The edges of $\cD$ belong neither to
$E(\cG)$ nor to $E(\cG^*)$. A star of a vertex $x_0\in V(\cG)$
produces a flower of adjacent quadrilaterals from $F(\cD)$ around
the common vertex $x_0$, see Fig.\,\ref{flower}.

Observe that the double $\cD$ is automatically {\it bipartite},
since its vertices $V(\cD)$ are decomposed into two complementary
halves, $V(\cD)=V(\cG)\sqcup V(\cG^*)$ (``black'' and ``white''
vertices), such that the ends of each edge from $E(\cD)$ are of
different colours. An arbitrary quad-graph embedded in $\bbC$ is
automatically bipartite, and the above construction can be
reversed for it, to produce a cell decomposition $\cG$ along with
its dual $\cG^*$. The decomposition of $V(\cD)$ into $V(\cG)$ and
$V(\cG^*)$ is unique, up to interchanging the roles of $\cG$ and
$\cG^*$. Edges of $\cG$ (say) connect two ``black'' vertices along
the diagonal of each face of $\cD$.

%-----------------------------------------------------------------
\begin{figure}[htbp]
\begin{minipage}[t]{200pt}
\setlength{\unitlength}{0.04em}
\begin{picture}(200,240)(-220,-140)
 \put(-100,0){\circle*{13}}\put(100,0){\circle*{13}}
 \put(0,-100){\circle{13}} \put(0,100){\circle{13}}
 \path(-100,0)(100,0)
 \path(-100,0)(-120,70)
 \path(-100,0)(-120,-70)
 \path(100,0)(120,70)
 \path(100,0)(120,-70)
 \dashline[+30]{10}(0,100)(0,-100)
 \dashline[+30]{10}(0,100)(-140, 30)
 \dashline[+30]{10}(0,100)(140, 30)
 \dashline[+30]{10}(0,-100)(-140, -30)
 \dashline[+30]{10}(0,-100)(140, -30)
 \thicklines
 \path(-100,0)(0,-100)
 \path(0,-100)(100,0)
 \path(100,0)(0,100)
 \path(0,100)(-100,0)
 \put(-140,  -7){$x_0$} \put(115, -7){$x_1$}
 \put(  -7,-125){$y_0$} \put( -7,115){$y_1$}
\end{picture}
\caption{A face of the double $\cD$}\label{diamond}
\end{minipage}\hfill
\begin{minipage}[t]{200pt}
    \setlength{\unitlength}{15pt}
\begin{picture}(7,9)(-2,0)
\put(4,4){\circle*{0.3}} \put(4.2,8){\circle*{0.3}}
\put(5,6){\circle{0.3}} \put(3,5.7){\circle{0.3}}
\put(7,5.9){\circle*{0.3}} \put(6,4.2){\circle{0.3}}
\put(7.5,2.6){\circle*{0.3}} \put(5.1,2.5){\circle{0.3}}
\put(4,0.2){\circle*{0.3}} \put(3,2.1){\circle{0.3}}
\put(0.4,3.9){\circle*{0.3}} \thicklines \path(4,4)(5,6)
\path(5,6)(4.2,8) \path(4.2,8)(3,5.7) \path(3,5.7)(4,4)
\path(3,5.7)(0.4,3.9) \path(0.4,3.9)(3,2.1) \path(3,2.1)(4,4)
\path(3,2.1)(4,0.2) \path(4,0.2)(5.1,2.5) \path(5.1,2.5)(4,4)
\path(5.1,2.5)(7.5,2.6) \path(7.5,2.6)(6,4.2) \path(6,4.2)(4,4)
\path(6,4.2)(7,5.9) \path(7,5.9)(5,6) \put(3.1,4){$x_0$}
\put(6.3,4.2){$y_0$}  \put(7.3,5.9){$x_1$}  \put(5.1,6.3){$y_1$}
\put(4.1,8.3){$x_2$}  \put(2.3,6){$y_2$}    \put(0,4.3){$x_3$}
\put(2.4,1.6){$y_3$}  \put(3.9,-0.3){$x_4$} \put(5.1,2.0){$y_4$}
\put(7.5,2.1){$x_5$}
\end{picture}
   \caption{Faces of $\cD$ around vertex $x_0$.}
        \label{flower}
         \end{minipage}
    \end{figure}
%------------------------------------------------------------------------

Let there be given a function $\nu:E(\cG)\to\bbC$ on undirected
edges of $\cG$. (It is assumed that both directed representatives
$\pm\ee$ of any edge carry the same value $\nu(\ee)=\nu(-\ee)$ as
the underlying undirected one.) Extend the function $\nu$ to
undirected edges of $\cG^*$ according to the rule
\begin{equation}\label{nu*}
\nu(\ee^*)=1/\nu(\ee).
\end{equation}
\begin{definition}
A function $f:V(\cD)\to\bbC$ is called {\itbf discrete
holomorphic} (with respect to the weights $\nu$), if for any
positively oriented quadrilateral $(x_0,y_0,x_1,y_1)\in F(\cD)$
there holds:
\begin{equation}\label{discr CR}
  \frac{f(y_1)-f(y_0)}{f(x_1)-f(x_0)}=i\nu(x_0,x_1)=-\frac{1}{i\nu(y_0,y_1)}.
\end{equation}
These equations are called {\itbf discrete Cauchy-Riemann
equations}.
\end{definition}
Again, the most interesting case corresponds to the real positive
weights $\nu:E(\cG)\sqcup E(\cG^*)\to\bbR_+$. The theory of
discrete holomorphic functions was developed in \cite{M1}. In
\cite{K} a discrete Dirac operator was introduced, the kernel
of which consists of discrete holomorphic functions. In the present
paper, like in \cite{M1}, our attention belongs not to the
discrete Dirac operator but to its kernel only. The next statement
follows immediately.
\begin{proposition}
a) If a function $f:V(\cD)\to\bbC$ is discrete holomorphic, then
its restrictions to $V(\cG)$ and to $V(\cG^*)$ are discrete
harmonic.

b) Conversely, any discrete harmonic function $f:V(\cG)\to\bbC$
admits a family of discrete holomorphic extensions to $V(\cD)$,
differing by an additive constant on $V(\cG^*)$. Such an extension
is uniquely defined by a value at one arbitrary vertex $y\in
V(\cG^*)$.
\end{proposition}

\section{Rhombic embeddings and labelings of a quad-graph}
\label{Sect KS}

The paper \cite{KS} studies {\it rhombic embeddings} of a
quad-graph $\cD$ in $\bbC$, i.e., embeddings with the property
that each face of $\cD$ is mapped to a rhombus. A combinatorial criterion
for the existence of a rhombic embedding of a given quad-graph
$\cD$ found in \cite{KS} is as follows.
\begin{definition}
A {\itbf strip} $\cS$ in $\cD$ is a path
$\{\aaa_j^*\}_{j=-\infty}^\infty$ in $\cD^*$ with the following
property: for any two consecutive dual edges
$\aaa_j^*,\aaa_{j+1}^*\in\cS\subset E(\cD^*)$ with the common point
$\aaa_j^*\cap\aaa_{j+1}^*=q_j\in V(\cD^*)\simeq F(\cD)$, the
corresponding edges $\aaa_j,\aaa_{j+1}\in E(\cD)$ are two
opposite sides of the quadrilateral $q_j$. The edges
$\{\aaa_j\}_{j=-\infty}^\infty$ are called the {\itbf traverse
edges} of the strip $\cS$.
\end{definition}
\begin{theorem}\label{theorem KS}
 {\rm \cite{KS}} A planar quad-graph $\cD$ admits a
rhombic embedding in $\bbC$ if and only if the following two
conditions are satisfied:
\begin{itemize}
\item No strip crosses itself or is periodic.
\item Two distinct strips cross each other at most once.
\end{itemize}
\end{theorem}
A rhombic embedding determines rhombus angles that are naturally
assigned to the edges of $\cG$ and $\cG^*$, see
Fig.\,\ref{rhombus}. Such systems of rhombus angles
$\phi:E(\cG)\sqcup E(\cG^*)\to (0,\pi)$ are characterized,
according to \cite{KS}, by the following two conditions: first,
\begin{equation}\label{phi*}
\phi(\ee^*)=\pi-\phi(\ee),\qquad \forall \ee\in E(\cG),
\end{equation}
and, second,
\begin{equation}\label{sum phi}
\sum_{\ee\in{\rm star}(x_0;\cG)}\phi(\ee)=2\pi, \quad
\sum_{\ee^*\in{\rm star}(y_0;\cG^*)}\phi(\ee^*)=2\pi,\qquad
\forall x_0\in V(\cG), \; y_0\in V(\cG^*).
\end{equation}

 %-----------------------------------------------------------------
\begin{figure}[htbp]
\begin{center}
\setlength{\unitlength}{0.04em}
\begin{picture}(200,220)(-100,-110)
 \put(-100,0){\circle*{13}}\put(100,0){\circle*{13}}
 \put(0,-75){\circle{13}} \put(0,75){\circle{13}}
 \put(-100,0){\line(4,-3){95}}
 \put(-100,0){\line(4,3){95}}
 \put(4,-72){\line(4,3){91}}
 \put(4,72){\line(4,-3){91}}
 \put(-137,  -7){$x_0$} \put(115, -7){$x_1$}
 \put(  -6,-100){$y_0$} \put( -6,97){$y_1$}
% \put(-75,57){$\alpha_1$} \put(60,-60){$\alpha_1$}
% \put(57,55){$\alpha_0$} \put(-75,-64){$\alpha_0$}
 \bezier{100}(-84,12)(-74,0))(-84,-12)
 \bezier{100}(-16,63)(0,53))(16,63)
 \bezier{100}(-12,66)(0,57))(12,66)
 \put(-68,  -7){$\phi(\ee)$}
 \put(-22,  35){$\phi(\ee^*)$}
\end{picture}
  \caption{A rhombic embedding of a quadrilateral
    $(x_0,y_0,x_1,z_1)\in F(\cD)$, $\ee=(x_0,x_1)$, $\ee^*=(y_0,y_1)$}
  \label{rhombus}
\end{center}
\end{figure}
%-----------------------------------------------------------------

As mentioned in \cite{KS}, to each rhombic embedding of $\cD$
there corresponds a set of parallelogram embeddings (wherein each
face is mapped to a parallelogram), which are obtained by
replacing each traverse edge of a strip with a real multiple (a
different multiple for each strip).

\begin{definition}\label{labelling}
A {\itbf labeling} is a function $\alpha:\vec{E}(\cD)\to\bbC$ such
that $\alpha(-\aaa)=-\alpha(\aaa)$ for any $\aaa\in\vec{E}(\cD)$, and
the values of $\alpha$ on two opposite and equally directed edges
of any quadrilateral from $F(\cD)$ are equal to one another.
\end{definition}
%-----------------------------------------------------------------
\begin{figure}[htbp]
\setlength{\unitlength}{0.04em}
\begin{minipage}[t]{200pt}
\begin{picture}(200,200)(-200,-110)
 \put(-100,0){\circle*{13}}\put(100,0){\circle*{13}}
 \put(0,-75){\circle{13}} \put(0,75){\circle{13}}
 \put(-96,-3){\vector(4,-3){91}}
 \put(-96,3){\vector(4,3){91}}
 \put(96,3){\vector(-4,3){91}}
 \put(96,-3){\vector(-4,-3){91}}
 \put(-137,  -7){$x_0$} \put(115, -7){$x_1$}
 \put(  -6,-100){$y_0$} \put( -6,97){$y_1$}
 \put(-75,55){$\alpha_1$} \put(50,-60){$-\alpha_1$}
 \put(45,55){$-\alpha_0$} \put(-75,-60){$\alpha_0$}
\end{picture}
\caption{Labeling of directed edges}\label{diamond again or}
\end{minipage}
\hfill
\begin{minipage}[t]{200pt}
\begin{picture}(200,200)(-200,-110)
 \put(-100,0){\circle*{13}}\put(100,0){\circle*{13}}
 \put(0,-75){\circle{13}} \put(0,75){\circle{13}}
 \put(-96,-3){\line(4,-3){91}}
 \put(-96,3){\line(4,3){91}}
 \put(4,-72){\line(4,3){91}}
 \put(4,72){\line(4,-3){91}}
 \put(-137,  -7){$x_0$} \put(115, -7){$x_1$}
 \put(  -6,-100){$y_0$} \put( -6,97){$y_1$}
 \put(-75,57){$\alpha_1^2$} \put(60,-60){$\alpha_1^2$}
 \put(57,55){$\alpha_0^2$} \put(-75,-64){$\alpha_0^2$}
\end{picture}
\caption{Labeling of undirected edges}\label{diamond again}
\end{minipage}
\end{figure}
%-----------------------------------------------------------------

This definition is illustrated in Fig.\,\ref{diamond again or}.
Note that if edges of any given face of $\cD$ are directed as on
this figure (from black to white), then any two opposite edges
carry opposite labels. For any labeling
$\alpha:\vec{E}(\cD)\to\bbC$ of directed edges, the function
$\alpha^2$ can be considered as a {\it labeling of undirected
edges}, i.e., as a function $\alpha^2:E(\cD)\to\bbC$ such that its
values on two opposite edges of any quadrilateral from $F(\cD)$
are equal to one another, see Fig.\,\ref{diamond again}.
Conversely, any labeling of undirected edges comes as a square of
some labeling of directed edges.

The existence of a labeling $\alpha:\vec{E}(\cD)\to\bbC$ is
equivalent to the existence of a {\it parallelogram immersion} of
the quad-graph $\cD$. Indeed, given a parallelogram immersion
$p:V(\cD)\to\bbC$, one defines canonically a labeling by
\begin{equation}\label{label canonical}
\alpha(x,y)=p(y)-p(x),\quad \forall (x,y)\in\vec{E}(\cD).
\end{equation}
Conversely, given a labeling $\alpha:\vec{E}(\cD)\to\bbC$, the
formula (\ref{label canonical}) correctly defines a function
$p:V(\cD)\to\bbC$ and assures that the $p$-image of any
quadrilateral face of $\cD$ is a parallelogram. If the labels
$\alpha$ take values in $\bbS^1=\{\theta\in{\Bbb C}:
|\theta|=1\}$, then the corresponding immersion is rhombic.

\begin{definition}
A parallelogram immersion $p:V(\cD)\to\bbC$ of a quad-graph $\cD$
is called {\itbf quasicrystallic}, if the set of values of the
corresponding labeling $\alpha:\vec{E}(\cD)\to\bbC$, defined by
(\ref{label canonical}), is finite, say
$A=\{\pm\alpha_1,\ldots,\pm\alpha_d\}$.
\end{definition}
It will be supposed that any two non-opposite elements of $A$ are
linearly independent over $\bbR$. This implies, in particular,
that all parallelograms are non-degenerate.

It will be of a central importance for us that any quasicrystallic
parallelogram immersion $p$ can be seen as a sort of a projection
of a certain two-dimensional subcomplex (combinatorial surface)
$\Omega_\cD$ of a multi-dimensional regular square lattice
$\bbZ^d$. The vertices of $\Omega_\cD$ are given by a map
$P:V(\cD)\to\bbZ^d$ constructed as follows. Fix some $x_0\in
V(\cD)$, and set $P(x_0)=\bO$. For all other vertices of $\cD$,
their images in $\bbZ^d$ are defined recurrently by the property:
\begin{itemize}
\item For any two neighbors $x,y\in V(\cD)$, if $p(y)-p(x)=\pm\alpha_i\in A$,
then $P(y)-P(x)=\pm\be_i\,$,
\end{itemize}
where $\be_i$ is the $i$-th coordinate vector of $\bbZ^d$. The
edges and faces of $\Omega_\cD$ correspond to edges and faces of
$\cD$. So, the combinatorics of $\Omega_\cD$ is that of $\cD$, and
therefore Theorem \ref{theorem KS} can be used to decide whether a
given two-dimensional subcomplex of $\bbZ^d$ corresponds in this
way to some rhombic embedding of a quad-graph in $\bbC$.

\section{3D consistency}
\label{Sect 3D}

We now study a question about {\em integrability} of the discrete
Cauchy-Riemann equations (\ref{discr CR}). These equations are
just a specific linear issue of general equations on quad-graphs
\cite{BS}
\begin{equation}\label{eq}
\Phi(f(x_0),f(y_0),f(x_1),f(y_1))=0,
\end{equation}
relating four fields $f$ sitting on the four vertices of an
arbitrary (oriented) face $(x_0,y_0,x_1,y_1)\in F(\cD)$ of a
quad-graph $\cD$. Here the function $\Phi$ may depend on some
parameters (in the case of discrete Cauchy-Riemann equations these
are the weights $\nu$), and it is supposed that equation
(\ref{eq}) is uniquely solvable for any one of the fields in terms
of other three (which is, of course, the case for discrete
Cauchy-Riemann equations with non-vanishing weights $\nu$).

The approach pushed forward in \cite{BS} is based on the idea
that integrability of such equations on quad-graphs is synonymous
with their 3D consistency. To describe the latter notion, we
extend the planar quad--graph $\cD$ into the third dimension.
Formally speaking, we consider the second copy $\widehat\cD$ of
$\cD$ and add edges connecting each vertex $x\in V(\cD)$ with its
copy $\widehat x\in V(\widehat\cD)$. On this way we obtain a
``three-dimensional quad-graph'' ${\bD}$, whose set of vertices is
\[
V({\bD})=V(\cD)\sqcup V(\widehat\cD),
\]
whose set of edges is
\[
E({\bD})=E(\cD)\sqcup E(\widehat\cD)\sqcup\{(x,\widehat x):\,x\in
V(\cD)\},
\]
and whose set of faces is
\[
F({\bD})=F(\cD)\sqcup F(\widehat\cD)\sqcup\{(x,y,\widehat
y,\widehat x): \,(x,y)\in E(\cD)\}.
\]
Elementary building blocks of $\bD$ are cubes
$(x_0,y_0,x_1,y_1,\widehat x_0,\widehat y_0,\widehat x_1,\widehat
y_1)$, as shown on Fig.\,\ref{cube}.
%-----------------------------------------------------------------
\begin{figure}[htbp]
\begin{center}
\setlength{\unitlength}{0.04em}
\begin{picture}(200,230)(0,0)
 \put(0,0){\circle*{10}}    \put(150,0){\circle{10}}
 \put(0,150){\circle{10}}   \put(150,150){\circle*{10}}
 \put(50,200){\circle*{10}} \put(200,200){\circle{10}}
 \put(50,50){\circle{10}}   \put(200,50){\circle*{10}}
 \path(0,0)(150,0)       \path(0,0)(0,150)
 \path(150,0)(150,150)   \path(0,150)(150,150)
 \path(0,150)(50,200)    \path(150,150)(200,200)   \path(50,200)(200,200)
 \path(200,200)(200,50) \path(200,50)(150,0)
 \dashline[+30]{10}(0,0)(50,50)
 \dashline[+30]{10}(50,50)(50,200)
 \dashline[+30]{10}(50,50)(200,50)
 \put(-33,-5){$x_0$} \put(-33,145){$\widehat x_0$}
 \put(160,-5){$y_0$} \put(160,140){$\widehat y_0$}
 \put(210,45){$x_1$} \put(210,200){$\widehat x_1$}
 \put(20,50){$y_1$}  \put(20,205){$\widehat y_1$}
\end{picture}
\caption{Elementary cube of $\bD$}\label{cube}
\end{center}
\end{figure}
%-----------------------------------------------------------------

Clearly, if $\cD$ is bipartite, then so is $\bD$: each $\widehat
x\in V(\widehat\cD)$ has the colour opposite to the colour of its
counterpart $x\in V(\cD)$. Hence, we can extend the ``black''
graph $\cG$ to a 3D object $\bG$, with edges
\[
E(\bG)=E(\cG)\sqcup E(\widehat{\cG^*})\sqcup \{(x,\widehat y):
x\in V(\cG),\;y\in V(\cG^*),\; (x,y)\in E(\cD)\}.
\]
Edges of $E(\bG)$ within Fig.\,\ref{cube} are $(x_0,x_1)$,
$(\widehat y_0,\widehat y_1)$, $(x_0,\widehat y_0)$,
$(x_0,\widehat y_1)$, $(x_1,\widehat y_0)$, and $(x_1,\widehat
y_1)$, forming the black tetrahedron. Similarly, we have a 3D
white graph $\bG^*$.
\begin{definition}
Equation (\ref{eq}) is called {\itbf 3D consistent} if it can be
imposed on all faces of any elementary cube of $\bD$, in such a
manner that opposite faces carry one and the same equation (i.e.,
the same parameters).
\end{definition}
This should be understood as follows. Consider an elementary cube
of $\bf D$, as on Fig.\,\ref{cube}. Suppose that the values of the
function $f$ are given at the vertex $x_0$ and at its three
neighbors $y_0$, $y_1$, and $\widehat x_0$. Then equation
(\ref{eq}) uniquely determines the values of $f$ at $x_1$,
$\widehat y_0$, and $\widehat y_1$. After that equation (\ref{eq})
delivers three {\it a priori} different values for the value of
the field $f$ at the vertex $\widehat x_1$, coming from the faces
$(y_0,x_1,\widehat x_1,\widehat y_0)$, $(x_1,y_1,\widehat
y_1,\widehat x_1)$, and $(\widehat x_0,\widehat y_0,\widehat
x_1,\widehat y_1)$, respectively. The 3D consistency means that
these three values for $f(\widehat x_1)$ actually coincide,
independently on the choice of initial conditions.

As discussed in detail in \cite{BS}, the 3D consistency of a given
system (\ref{eq}) allows one to construct B\"acklund
transformations and to find in an algorithmic way a zero curvature
representation for it, which are traditionally considered as main
attributes of integrability. Briefly, the constructions are as
follows.

1) Given a solution $f:V(\cD)\to\bbC$ to (\ref{eq}) and an
arbitrary value $f(\widehat x_0)=\widehat f_0$ at some vertex
$\widehat x_0\in\widehat\cD$, the 3D consistency allows one to
extend the solution $f$ to the whole of $V(\bD)$. Its restriction
to $V(\widehat\cD)$ is thus a well-defined function
$f:V(\widehat\cD)\to\bbC$ which also solves the original equation
(\ref{eq}). Setting $\widehat f(x)=f(\widehat x)$ for all $x\in
V(\cD)$, one can interpret this function as $\widehat
f:V(\cD)\to\bbC$, and this $\widehat f$ is called the {\itbf
B\"acklund transformation} of $f$ (defined by the value $\widehat
f_0$ and the parameters sitting on the vertical faces).

2) Suppose that the function $\Phi(u_1,u_2,u_3,u_4)$ in (\ref{eq})
is affine-linear in all its arguments, so that this equation can
be solved uniquely for an arbitrary argument $u_i$ in terms of
other three arguments, the solution being given by a
fractional-linear function. For an arbitrary edge
$\aaa=(x,y)\in\vec{E}(\cD)$, consider the vertical face
$(x,y,\widehat y, \widehat x)\in F(\bD)$ over this edge. The
solution of equation $\Phi(f(x),f(y),f(\widehat y),f(\widehat
x))=0$ can be written as
\begin{equation}\label{L prelim}
f(\widehat y)=L(f(y),f(x))\cdot f(\widehat x),
\end{equation}
where $L(f(x),f(y))\in PGL_2(\bbC)$, and the standard notation for
the action of $PGL_2(\bbC)$ on $\bbC$ by M\"obius transformations
is used:
\[
\begin{pmatrix} a & b \\ c & d \end{pmatrix}\cdot
u=\dfrac{au+b}{cu+d}\,.
\]
One assigns the matrix above to the edge $\aaa$, so that
$L(\aaa)=L(f(y),f(x))$. Now it follows from the 3D consistency that
for an arbitrary face $(x_0,y_0,x_1,y_1)\in F(\cD)$ one has:
\begin{equation}\label{zcr prelim}
L(f(x_1),f(y_0))\,L(f(y_0),f(x_0))=L(f(x_1),f(y_1))\,L(f(y_1),f(x_0)).
\end{equation}
This expresses the flatness of the discrete connection $L$ on
$\cD$ with values in $PGL_2(\bbC)$, hence (\ref{zcr prelim}) is
called the {\itbf zero  curvature representation} of system
(\ref{eq}). It is often possible to use suitable normalizations in
order to lift this representation to the one with values in
$GL_2(\bbC)$.

\section{3D consistent Cauchy-Riemann equations}
\label{Sect 3D CR}

To apply the notion of the 3D consistency to the discrete
Cauchy-Riemann equation (\ref{discr CR}), one has to explain how
to impose it on the vertical faces of $\bD$. For this, we assume
that the function $\nu$ is extended to $E(\bG)\sqcup E(\bG^*)$,
still satisfying the condition $\nu(\ee^*)=1/\nu(\ee)$, and with
an additional condition that {\em opposite edges carry the same
values of $\nu$}.

An interesting problem is, of course, to find functions $\nu$ on
the ``ground floor'' $E(\cG)\sqcup E(\cG^*)$ which can be extended
to the edges of $E(\bG)\sqcup E(\bG^*)$ lying in the ``vertical''
faces to give a 3D consistent system.
\begin{theorem}\label{Th 3D CR}
The function $\nu:E(\cG)\sqcup E(\cG^*)\to\bbC$ can be extended to
$E(\bG)\sqcup E(\bG^*)$ giving a 3D consistent system of discrete
Cauchy-Riemann equations, if and only if the following condition
is satisfied:
\begin{equation}\label{prod=1}
\prod_{\ee\in{\rm star}(x_0;\cG)}
 \frac{1+i\nu(\ee)}{1-i\nu(\ee)}=1,\quad
 \prod_{\ee^*\in{\rm star}(y_0;\cG^*)}
 \frac{1+i\nu(\ee^*)}{1-i\nu(\ee^*)}=1,\qquad
 \forall x_0\in V(\cG), \; y_0\in V(\cG^*).
\end{equation}
\end{theorem}
{\bf Proof.} Consider a flower of quadrilaterals around $x_0$,
with $\ee_k=(x_0,x_k)$, $\ee_k^*=(y_{k-1},y_k)$ (in notations of
Fig.\,\ref{flower}). Build the extension of this flower to the
third dimension (for one of its petals corresponding to $k=1$ this
extension is shown on Fig.\,\ref{cube}). Denote
\begin{equation}\label{nu's}
  \nu(y_{k-1},y_k)=\nu(\ee_k^*)=\nu_k, \quad \nu(y_k,\widehat x_0)=\mu_k.
\end{equation}
\begin{lemma}\label{lemma 3D nus}
Discrete Cauchy-Riemann equations are 3D consistent on the cube
over the $k$-th petal, if and only if
\begin{equation}\label{3D CR aux}
  1+\nu_k\mu_{k-1}-\nu_k\mu_k+\mu_{k-1}\mu_k=0.
\end{equation}
\end{lemma}
{\bf Proof of Lemma \ref{lemma 3D nus}.} Consider the elementary
cube on Fig.\,\ref{cube}, corresponding to $k=1$. On the first
step of checking the 3D consistency we find:
\begin{eqnarray*}
 f(x_1)  & = & f(x_0)+i\nu_1(f(y_0)-f(y_1)),\\
 f(\widehat y_0) & = & f(x_0)+i\mu_0(f(y_0)-f(\widehat x_0)),\\
 f(\widehat y_1) & = & f(x_0)+i\mu_1(f(y_1)-f(\widehat x_0)).
\end{eqnarray*}
On the second step we find (from the condition that opposite faces
support the same equations):
\begin{eqnarray*}
 f(\widehat x_1)  & = &
                  f(\widehat x_0)+i\nu_1(f(\widehat y_0)-f(\widehat y_1)),\\
          & = & f(y_1)+i\mu_0(f(x_1)-f(\widehat y_1)),\\
          & = & f(y_0)+i\mu_1(f(x_1)-f(\widehat y_0)).
\end{eqnarray*}
After simple computations we find:
\begin{eqnarray*}
 f(\widehat x_1)  & = &
  (1+\nu_1\mu_0-\nu_1\mu_1)f(\widehat x_0)
  -\nu_1\mu_0f(y_0)+\nu_1\mu_1f(y_1),\\
          & = &
  -\mu_0\mu_1f(\widehat x_0)-\nu_1\mu_0f(y_0)
  +(1+\nu_1\mu_0+\mu_0\mu_1)f(y_1),\\
          & = &
  -\mu_0\mu_1f(\widehat x_0)+(1-\nu_1\mu_1+\mu_0\mu_1)f(y_0)
  +\nu_1\mu_1f(y_1).
\end{eqnarray*}
Comparison of these expressions leads to
$1+\nu_1\mu_0-\nu_1\mu_1+\mu_0\mu_1=0$, which proves the lemma.
\qed
\medskip

Continuing the proof of Theorem \ref{Th 3D CR}, we derive from
(\ref{3D CR aux}):
\[
\mu_k=\frac{\nu_k\mu_{k-1}+1}{\nu_k-\mu_{k-1}}=\begin{pmatrix}\nu_k
& 1 \\ -1 & \nu_k \end{pmatrix}\cdot \mu_{k-1},
\]
where the standard notation for the action of $PGL_2(\bbC)$ on
$\bbC$ by M\"obius transformations is used. Starting with an
arbitrary $\mu_0$, we can define all $\mu_k$'s consecutively. This
procedure is consistent, if running around $x_0$ returns the value
of $\mu_0$ we started with. This holds for any $\mu_0$, if and
only if the matrix product
$
\displaystyle{\prod_k^{\curvearrowleft}}\begin{pmatrix}\nu_k & 1
\\ -1 & \nu_k \end{pmatrix}
$
is a scalar matrix. It is easy to see by induction that the above
matrix product may be presented as
\[
\begin{pmatrix}A & B \\ -B & A \end{pmatrix}\quad {\rm with}\quad
A=\frac{1}{2}\Big(\prod_k(\nu_k+i)+\prod_k(\nu_k-i)\Big),\quad
B=\frac{1}{2i}\Big(\prod_k(\nu_k+i)-\prod_k(\nu_k-i)\Big).
\]
Therefore, a necessary and sufficient condition for this matrix to
be scalar is
\[
B=0\quad\Leftrightarrow\quad\prod_k\frac{\nu_k+i}{\nu_k-i}=1,
\]
which is equivalent to the first equality in (\ref{prod=1}),
because of $\nu_k=\nu(\ee_k^*)=1/\nu(\ee_k)$. The second condition
in (\ref{prod=1}) is proved similarly, by considering a flower of
quadrilaterals around $y_0\in V(\cG^*)$. \qed
\medskip

As pointed out above, the most interesting case is when $\nu$
takes values in $\bbR_+$. In this case we will use the notation
\begin{equation}\label{phi}
\nu(\ee)=\tan\frac{\phi(\ee)}{2},\quad \phi(\ee)\in(0,\pi).
\end{equation}
The condition $\nu(\ee^*)=1/\nu(\ee)$ is translated in this case
into (\ref{phi*}). The integrability condition (\ref{prod=1})
takes in this case the form
\begin{equation}\label{prodexp=1}
\prod_{\ee\in{\rm star}(x_0;\cG)}\exp(i\phi(\ee))=1, \quad
\prod_{\ee^*\in{\rm star}(y_0;\cG^*)}\exp(i\phi(\ee^*))=1,\qquad
\forall x_0\in V(\cG),\; y_0\in V(\cG^*).
\end{equation}
The latter condition is a generalization of (\ref{sum phi}), and
is equivalent to saying that the system of angles
$\phi:E(\cG)\sqcup E(\cG^*)\to(0,\pi)$ comes from a realization of
the quad-graph $\cD$ by a {\it rhombic ramified embedding} in $\bbC$.
Flowers of such an embedding can wind around its vertices more
than once.

\begin{lemma}\label{lemma on factorization or}
Let a quad-graph $\cD$ be a double for a pair of dual cell
decompositions $\cG$, $\cG^*$. Let $\Phi:E(\cG)\sqcup
E(\cG^*)\to\bbC$ be a function satisfying
\begin{equation}\label{Phi*}
\Phi(\ee^*)=-1/\Phi(\ee), \qquad\forall \ee\in E(\cG).
\end{equation}
Then the necessary and sufficient condition for the existence of a
labeling $\alpha:\vec{E}(\cD)\to\bbC$ such that, in the notations
of Fig.\,\ref{diamond again or},
\begin{equation}\label{lemma 0 or}
\Phi(\ee)=\Phi(x_0,x_1)=\frac{\alpha_1}{\alpha_0}\quad\Leftrightarrow\quad
\Phi(\ee^*)=\Phi(y_0,y_1)=-\frac{\alpha_0}{\alpha_1}\,,
\end{equation}
is given by the equations
\begin{equation}\label{lemma 1 or}
\prod_{\ee\in{\rm star}(x_0;\cG)} \Phi(\ee)=1,\quad
 \prod_{\ee^*\in{\rm star}(y_0;\cG^*)} \Phi(\ee^*)=1,\qquad
 \forall x_0\in V(\cG), \; y_0\in V(\cG^*).
\end{equation}
\end{lemma}
{\bf Proof.} The necessity is obvious. To prove sufficiency, we
construct $\alpha$ by assigning an arbitrary value (say,
$\alpha=1$) to some edge of $\cD$, and then extending it
successively using either of the equations (\ref{lemma 0 or}) and
the definition of labeling. Conditions (\ref{lemma 1 or}) assure
the consistency of this procedure. \qed

\begin{corollary}
Integrability condition (\ref{prod=1}) for the function
$\nu:E(\cG)\sqcup E(\cG^*)\to\bbC$ is equivalent to the following
one: there exists a labeling $\alpha:\vec{E}(\cD)\to\bbC$ of
directed edges of $\cD$, such that, in notations of
Fig.\,\ref{diamond again or},
\begin{equation}\label{nu paral}
\nu(y_0,y_1)=\frac{1}{\nu(x_0,x_1)}=
i\,\frac{\alpha_1+\alpha_0}{\alpha_1-\alpha_0}.
\end{equation}
Under this condition, the 3D consistency of the discrete
Cauchy-Riemann equations is assured by the following values of the
weights $\nu$ on the edges of $E(\bG)\sqcup E(\bG^*)$ lying in the
vertical faces:
\begin{equation}\label{nu vert}
\nu(y,\widehat{x})=\frac{1}{\nu(x,\widehat{y})}=
i\,\frac{\lambda+\alpha}{\lambda-\alpha},
\end{equation}
where $\alpha=\alpha(x,y)$, and $\lambda\in\bbC$ is an arbitrary
number having the interpretation of the label carried by all
vertical edges of $\bD$:
$\lambda=\alpha(x,\widehat{x})=\alpha(y,\widehat{y})$.
\end{corollary}
{\bf Proof.} Apply Lemma \ref{lemma on factorization or} with the
function
\begin{equation}\label{Phi CR}
\Phi(\ee)=\frac{1+i\nu(\ee)}{1-i\nu(\ee)},
\end{equation}
which satisfies (\ref{Phi*}) due to the property (\ref{nu*}) of
the weights $\nu$. Note that in the case $\nu(\ee)\in\bbR_+$ the
notation (\ref{phi}) implies that $\Phi(\ee)=\exp(i\phi(\ee))$.
The formula (\ref{lemma 0 or}) with the function (\ref{Phi CR}) is
clearly equivalent to (\ref{nu paral}). To prove the second
statement, we use notations of Lemma \ref{lemma 3D nus}, in
particular the formula $1+\nu_1\mu_0-\nu_1\mu_1+\mu_0\mu_1=0$.
According to (\ref{nu paral}), we have:
$\nu_1=i\frac{\alpha_1+\alpha_0}{\alpha_1-\alpha_0}$. Parametrize
the (arbitrary) value of $\mu_0$ as
$\mu_0=i\frac{\lambda+\alpha_0}{\lambda-\alpha_0}$. Then it
follows from the above formula that
$\mu_1=i\frac{\lambda+\alpha_1}{\lambda-\alpha_1}$. An easy
induction proves (\ref{nu vert}) for all edges in the vertical
faces.  \qed
\medskip

So, integrability of the discrete Cauchy-Riemann equations is
equivalent to the existence of a labeling $\alpha$ of directed
edges satisfying (\ref{nu paral}). Let $p:V(\cD)\to\bbC$ be a
parallelogram realization of $\cD$ defined by
$p(y)-p(x)=\alpha(x,y)$. Then discrete holomorphic functions are
characterized by
\begin{equation}\label{CR paral}
\frac{f(y_1)-f(y_0)}{f(x_1)-f(x_0)}=
\frac{\alpha_1-\alpha_0}{\alpha_1+\alpha_0}=
\frac{p(y_1)-p(y_0)}{p(x_1)-p(x_0)}.
\end{equation}
In other words, the quotient of diagonals of the $f$-image of any
quadrilateral $(x_0,y_0,x_1,y_1)\in F(\cD)$ is equal to the
quotient of diagonals of the corresponding parallelogram. In the
case of positive weights $\nu\in\bbR_+$, the labels $\alpha$ take
values in $\bbS^1$, and have a geometric interpretation of edges
of a rhombic realization of $\cD$.

\begin{proposition}\label{prop CR zcr}
The discrete Cauchy-Riemann equations (\ref{CR paral}) admit a
zero curvature representation (\ref{zcr prelim}) in
$GL_2(\bbC)[\lambda]$, with transition matrices along
$(x,y)\in\vec{E}(\cD)$ given by
\begin{equation}\label{L CR}
\renewcommand{\arraystretch}{1.4}
L(y,x,\alpha;\lambda)=\begin{pmatrix} \lambda+\alpha &
-2\alpha(f(x)+f(y)) \\ 0 & \lambda-\alpha \end{pmatrix}, \quad
where \quad\alpha=p(y)-p(x).
\end{equation}
\end{proposition}
{\bf Proof.} This result is easy to check. It can be also
systematically derived using the procedure outlined at the end of
Sect. \ref{Sect 3D}. Indeed, setting $\lambda=p(\widehat x)-p(x)$,
one writes the equation (\ref{CR paral}) on the vertical face
$(x,y,\widehat y,\widehat x)$ as
\[
\frac{f(\widehat x)-f(y)}{f(\widehat y)-f(x)}=
\frac{\lambda-\alpha}{\lambda+\alpha}\quad \Leftrightarrow \quad
f(\widehat y)= \frac{\lambda+\alpha}{\lambda-\alpha}\,f(\widehat
x)+\Big(f(x)- \frac{\lambda+\alpha}{\lambda-\alpha}\,f(y)\Big)=
M(y,x,\alpha;\lambda)\cdot f(\widehat x),
\]
where
\[
\renewcommand{\arraystretch}{1.4}
M(y,x,\alpha;\lambda)=\begin{pmatrix} \lambda+\alpha &
(\lambda-\alpha)f(x)-(\lambda+\alpha)f(y) \\ 0 & \lambda-\alpha
\end{pmatrix}.
\]
One easily shows that these matrices form a zero curvature
representation with values in $GL_2(\bbC)[\lambda]$, i.e., that
(\ref{zcr prelim}) holds literally, and not only projectively (up
to a scalar factor). Finally, observe that the matrices $L$ in
(\ref{L CR}) are gauge equivalent to the matrices $M$:
\[
L(y,x,\alpha;\lambda)=\begin{pmatrix} 1 & -f(y) \\ 0 & 1
\end{pmatrix}M(y,x,\alpha;\lambda)\begin{pmatrix} 1 & f(x) \\ 0 &
1 \end{pmatrix}.
\]
This finishes the proof. \qed
\medskip

The main result of the present section can be formulated as
follows. {\em Discrete Cauchy-Riemann equations on a quad-graph
$\cD$ are integrable if and only if they come from a parallelogram
immersion of $\cD$ in $\bbC$, weights $i\nu$ being the quotients
of diagonals of the corresponding parallelograms. In the case of
real positive weights $\nu$ the parallelograms are actually
rhombi.}

\section{Extension of discrete holomorphic functions to $\bbZ^d$}
\label{Sect CR multidim}

\begin{figure}[htbp]
\begin{center}
\setlength{\unitlength}{0.04em}
\begin{picture}(500,230)(-80,-50)
 \path(10,10)(-70,-50)(70,-50)(70,100)(280,100)(280,0)(350,0)(430,60)
      (360,60)(360,160)(150,160)(70,100)
 \path(110,130)(320,130)(320,30)(390,30)
 \path(140,100)(220,160)
 \path(210,100)(290,160)
 \path(280,100)(360,160)
 \path(280, 50)(360,110)
 \path(280,  0)(360, 60)
 \dashline[+30]{10}(150,160)(150,100)
 \path(150,100)(150,10)
 \path(70,-50)(150,10)
 \path(10,10)(70,10)
 \dashline[+30]{10}(70,10)(150,10)
 \path(-30,-20)(70,-20)
 \dashline[+35]{10}(70,-20)(110,-20)
 \path(110,-20)(110,100)
 \dashline[+30]{10}(110,100)(110,130)
 \path(0,-50)(70,2.5)
 \path(70,50)(136.6,100)
 \path(70, 0)(150, 60)
% \path(0,0)(150,0)       \path(0,0)(0,150)
% \path(150,0)(150,150)   \path(0,150)(150,150)
% \path(0,150)(50,200)    \path(150,150)(200,200)   \path(50,200)(200,200)
% \path(200,200)(200,50) \path(200,50)(150,0)
% \dashline[+30]{10}(0,0)(50,50)
% \dashline[+30]{10}(50,50)(50,200)
% \dashline[+30]{10}(50,50)(200,50)
\end{picture}
\caption{A non-monotone surface in $\bbZ^3$}\label{fig nonconvex}
\end{center}
\end{figure}
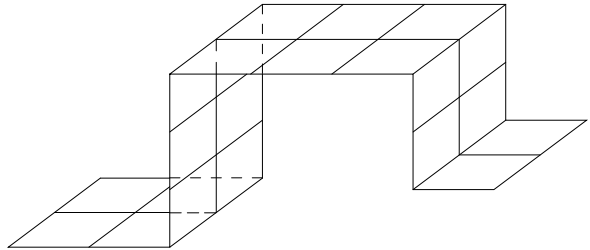

To exploit analytic possibilities provided by 3D consistency of
the discrete Cauchy-Riemann equations, we restrict our
considerations to quasicrystallic rhombic embeddings $\cD$, with
the set of labels $A=\{\pm\alpha_1,\ldots,\pm\alpha_d\}$.
Construct the two-dimensional subcomplex $\Omega_\cD$ in $\bbZ^d$
corresponding to $\cD$, as explained at the end of Sect. \ref{Sect KS}. 
Extend the labeling
$\alpha:\vec{E}(\cD)\to\bbC$ to all edges of $\bbZ^d$, assuming
that all edges parallel to (and directed as) $\be_k$ carry the
label $\alpha_k$. Now, 3D consistency of the discrete
Cauchy-Riemann equations allows us to impose them not only on
$\Omega_\cD$, but on the whole of $\bbZ^d$.
\begin{definition}\label{def multidim holo}
A function $f:\bbZ^d\to\bbC$ is called {\itbf discrete
holomorphic}, if it satisfies, on each elementary square of
$\,\bbZ^d$, the equation
\begin{equation}\label{CR square}
  \frac{f(\bn+\be_j+\be_k)-f(\bn)}{f(\bn+\be_j)-f(\bn+\be_k)}=
  \frac{\alpha_j+\alpha_k}{\alpha_j-\alpha_k}\,.
\end{equation}
\end{definition}
Obviously, for any discrete holomorphic function
$f:\bbZ^d\to\bbC$, its restriction to $V(\Omega_\cD)\sim V(\cD)$
is a discrete holomorphic function on $\cD$. To justify the
reverse procedure, i.e., the extension of an arbitrary discrete
holomorphic function on $\cD$ to $\bbZ^d$, keeping the property of
being discrete holomorphic, more thorough considerations are
necessary. 
\begin{definition}
For a given set $V\subset\bbZ^d$, its {\itbf hull} $\cH(V)$  is
the minimal set $\cH\subset\bbZ^d$ containing $V$ and satisfying
the condition: if three vertices of an elementary square belong to
$\cH$, then so does the fourth vertex.
\end{definition}
This notion is tailored for arbitrary 3D consistent four-point
equations of the type (\ref{eq}), including the discrete
Cauchy-Riemann equations. It is not difficult to show by induction
that the hull of an arbitrary connected subcomplex of $\bbZ^d$ is
a {\it brick}, i.e., a set of the type
\begin{equation}\label{brick}
 \Pi_{\ba,\bb}=\big\{\bn=(n_1,\ldots,n_d)\in\bbZ^d:\;
 a_k\le n_k\le b_k,\;\, k=1,\ldots,d\,\big\}\,,
\end{equation}
where $\ba=(a_1,\ldots,a_d)$, $\bb=(b_1,\ldots,b_d)$ are some
integer vectors, with infinite values $a_k=-\infty$, $b_k=\infty$
allowed. (Observe that Definition \ref{def multidim holo} is
equally well applicable to functions on bricks.) 
However, there exist combinatorial surfaces $\Omega$ (two-dimensional 
subcomplexes of $\bbZ^d$), like the one shown on Fig.\,\ref{fig nonconvex}, 
that support discrete holomorphic function which cannot be extended 
to $\bbZ^d$: the recursive process of extending an 
arbitrary discrete holomorphic function from $V(\Omega)$ to its hull
$\cH(V(\Omega))$ will lead to contradictions. The reason for this is
a non-monotonicity of $\Omega$: it contains pairs of points which cannot
be connected by a path in $\Omega$ with all edges lying in one octant.
However, such surfaces do not come 
from rhombic embeddings. We will prove 
the absence of contradictions in the case of $\Omega_\cD$.
%-----------------------------------------------------------------
%-----------------------------------------------------------------

\begin{proposition}\label{prop extend}
For a combinatorial surface $\Omega_\cD$ in $\bbZ^d$ coming from a 
rhombic embedding of a quad-graph $\cD$, set
\begin{equation}
 a_k=a_k(\Omega_\cD)=\min_{\bn\in V(\Omega_\cD)} n_k,\qquad
 b_k=b_k(\Omega_\cD)=\max_{\bn\in V(\Omega_\cD)} n_k,\qquad
 k=1,\ldots,d.
\end{equation}
(In case that $n_k$ are unbounded from below or from above on
$V(\Omega_\cD)$, set $a_k(\Omega_\cD)=-\infty$, resp.
$b_k(\Omega_\cD)=\infty$.) Then
$\cH(V(\Omega_\cD))=\Pi_{\ba,\bb}$, and an arbitrary discrete
holomorphic function on $\Omega_\cD$ can be extended to a discrete
holomorphic function on $\Pi_{\ba,\bb}$ in a unique and
unambiguous way.
\end{proposition}
For a proof of this proposition, a more detailed study of the
surface $\Omega_\cD$ will be necessary. In order to fix the ideas,
we will assume, without loss of generality, that the circular
order of the points $\pm\alpha_k$ on the positively oriented unit
circle $\bbS^1$ is the following:
$\alpha_1,\ldots,\alpha_d,-\alpha_1,\ldots,-\alpha_d$. We set
$\alpha_{k+d}=-\alpha_k$ for $k=1,\ldots,d$, and then define
$\alpha_m$ for all $m\in\bbZ$ by $2d$-periodicity.

Consider the set $A_m=\{\alpha_m,\ldots,\alpha_{m+d-1}\}$ of $d$
consecutive edge slopes. The opening angle of the sector spanned by
$\alpha_m$ and $\alpha_{m+d-1}$ is in $(0,\pi)$. The set $A_m$ contains 
exactly one member $\epsilon_k\alpha_k$ of each pair $\pm\alpha_k$, 
$k=1,\ldots,d$. This associates to any $m\in\bbZ$ the set of signs
$\beps=(\epsilon_1,\ldots,\epsilon_d)$, $\epsilon_k=\pm 1$, which
will be denoted by $\beps(m)$. The sets of signs $\beps(m)$ repeat
$2d$-periodically, therefore not all possible sets of signs appear
among them, but only the following $2d$ different ones. If
$m\in[1,d]$, then the corresponding $\beps=\beps(m)$ is given by
\[
\epsilon_k(m)=\left\{\begin{array}{ll} -1, & 1\le k<m,\\
 +1, & m\le k\le d,\end{array}\right.
\]
and if $m\in[d+1,2d]$, then
\[
\epsilon_k(m)=\left\{\begin{array}{ll} +1, & 1\le k<m-d, \\
 -1, & m-d\le k\le d.\end{array}\right.
\]

Fix an arbitrary $x_0\in V(\cD)$, and define the ``sector'' $U_m$
on the embedding plane $\bbC$ of the quad-graph $\cD$ as the set
of all points of $V(\cD)$ which can be reached from $x_0$ along
paths with all edges from
$A_m=\{\alpha_m,\ldots,\alpha_{m+d-1}\}$.

This can be re-formulated in terms of $\Omega_\cD$ as follows.
Recall that the map $P$ which identifies $\Omega_\cD$ with $\cD$
depends on the choice of the point $x_0\in V(\cD)$ corresponding
to $\bO\in V(\Omega_\cD)$. The map $P$ sends $U_m$ to the set of
the points of $V(\Omega_\cD)$ which can be reached from $\bO$
along paths in $V(\Omega_\cD)$ with all edges from
$\{\epsilon_1\be_1,\ldots,\epsilon_d\be_d\}$, where $\epsilon_k=
\epsilon_k(m)$ for $k=1,\ldots,d$. To formulate it in a still
another way, put into a correspondence to any set of signs
$\beps=(\epsilon_1,\ldots,\epsilon_d)$ the $d$-dimensional octant
\begin{equation}\label{sector}
  S_{\beps}=(\epsilon_1\bbZ_+)\times\ldots\times(\epsilon_d\bbZ_+)
  \subset\bbZ^d.
\end{equation}
In case of $\beps=\beps(m)$, use the notation $S_{\beps(m)}=S_m$.
Then the definition of $U_m$ is equivalent to saying that
$U_m=P^{-1}(V(\Omega_\cD)\cap S_m)$. The following statement will
be of a key importance.
\begin{lemma}\label{lemma sectors}
The union $\bigcup_{m=1}^{2d}U_m$ covers the whole of the
quad-graph $\cD$. Equivalently, the combinatorial surface
$\Omega_\cD$ coming from a rhombic embedding of $\cD$ lies
entirely in $\bigcup_{m=1}^{2d}S_m$.
\end{lemma}
{\bf Proof.} Clearly, $U_m$ lies within the sector of the embedding plane
with the tip at $x_0$, spanned by the directions $\alpha_m$ and 
$\alpha_{m+d-1}$. The set $A_m$ can be ordered:
$\alpha_m\prec\ldots\prec\alpha_{m+d-1}$. The lower boundary
$U_m^-$ (upper boundary $U_m^+$) of $U_m$ can be described as the
path in $\cD$ from the point $x_0$ obtained by following, at each
vertex of the path, the edge with the least (resp. the largest)
slope from $A_m$ available at this vertex, with respect to the
above mentioned ordering in $A_m$. The fact that $\cD$ is embedded
implies that all vertices of $\cD$ between $U_m^-$ and $U_m^+$
belong to $U_m$. Indeed, suppose that there are vertices between
$U_m^-$ and $U_m^+$ which cannot be reached from $x_0$ along a
path with all edges from $A_m$. Take such a vertex $x$,
combinatorially nearest to $x_0$. It cannot be reached from
$x_0$ along a path with the {\it last} edge from $A_m$. Then one
of the corners of one of the faces adjacent to $x$ is free from
edges from $A_m$ and therefore has an internal angle larger than
$\pi$, in a contradiction with embeddedness. Thus, $U_m$ can be
described as a set of vertices between $U_m^-$ and $U_m^+$.
Further, observe that the boundaries of the sectors $U_m$ are
interlaced: $U_m$ contains all $U_r^-$ with $r\in[m+1,m+d-2]$, and 
all $U_r^+$ with $r\in[m-d+2,m-1]$. This yields that the union of all 
$U_m$'s covers the whole of $\cD$. \qed
\smallskip

\begin{figure}[htb]
\begin{center}
\setlength{\unitlength}{0.15em}
%-----------------------------------
\begin{picture}(220,140)(-30,-20)
\multiput(0,0)(0,30){4}{
\multiput(0,0)(17.32,0){9}{\dottedline{2}(0,0)(8.66,5)
\dottedline{2}(0,0)(-8.66,5)\dottedline{2}(0,0)(0,-10)}}
\multiput(8.66,15)(17.32,0){9}{\dottedline{2}(0,0)(8.66,5)
\dottedline{2}(0,0)(-8.66,5)\dottedline{2}(0,0)(0,-10)}
\multiput(8.66,45)(17.32,0){9}{\dottedline{2}(0,0)(8.66,5)
\dottedline{2}(0,0)(-8.66,5)\dottedline{2}(0,0)(0,-10)}
\multiput(8.66,75)(17.32,0){9}{\dottedline{2}(0,0)(8.66,5)
\dottedline{2}(0,0)(-8.66,5)\dottedline{2}(0,0)(0,-10)}
\multiput(8.66,105)(17.32,0){9}{\dottedline{2}(0,0)(8.66,5)
\dottedline{2}(0,0)(-8.66,5)\dottedline{2}(0,0)(0,-10)}
\multiput(0,10)(17.32,0){9}{\dottedline{2}(0,0)(8.66,-5)
\dottedline{2}(0,0)(-8.66,-5)\dottedline{2}(0,0)(0,10)}
\multiput(0,40)(17.32,0){9}{\dottedline{2}(0,0)(8.66,-5)
\dottedline{2}(0,0)(-8.66,-5)\dottedline{2}(0,0)(0,10)}
\multiput(0,70)(17.32,0){9}{\dottedline{2}(0,0)(8.66,-5)
\dottedline{2}(0,0)(-8.66,-5)\dottedline{2}(0,0)(0,10)}
\multiput(0,100)(17.32,0){9}{\dottedline{2}(0,0)(8.66,-5)
\dottedline{2}(0,0)(-8.66,-5)\dottedline{2}(0,0)(0,10)}
\multiput(8.66,-5)(17.32,0){9}{\dottedline{2}(0,0)(8.66,-5)
\dottedline{2}(0,0)(-8.66,-5)\dottedline{2}(0,0)(0,10)}
\multiput(8.66,25)(17.32,0){9}{\dottedline{2}(0,0)(8.66,-5)
\dottedline{2}(0,0)(-8.66,-5)\dottedline{2}(0,0)(0,10)}
\multiput(8.66,55)(17.32,0){9}{\dottedline{2}(0,0)(8.66,-5)
\dottedline{2}(0,0)(-8.66,-5)\dottedline{2}(0,0)(0,10)}
\multiput(8.66,85)(17.32,0){9}{\dottedline{2}(0,0)(8.66,-5)
\dottedline{2}(0,0)(-8.66,-5)\dottedline{2}(0,0)(0,10)}
 \put(95.26,45){\circle*{3}}
% Sector S+++
 \thicklines
 \path(95.26,45)(86.6,50)
 {\thinlines
    \multiput(93.8744,45.8)(-1.3856,0.8){6}{\path(0,0)(1.73,1)}}
  \multiput(86.6,50)(-8.66,15){4}{\path(0,0)(-8.66,5)(-8.66,15)
  {\thinlines
    \multiput(-1.3856,0.8)(-1.3856,0.8){5}{\path(0,0)(1.73,1)}
    \multiput(-8.66,5)(0,1.6){5}{\path(0,0)(1.73,1)}}}
 \path(95.26,45)(103.92,50)
 {\thinlines
    \multiput(95.26,45)(1.3856,0.8){7}{\path(0,0)(-1.73,1)}}
  \multiput(103.92,50)(8.66,15){4}{\path(0,0)(8.66,5)(8.66,15)
  {\thinlines
    \multiput(1.0392,0.6)(1.3856,0.8){6}{\path(0,0)(-1.73,1)}
    \multiput(8.66,5.4)(0,1.6){7}{\path(0,0)(-1.73,1)}}}
% Sector S-++
  \path(95.26,45)(86.6,50)
   \multiput(86.6,50)(-8.66,15){4}{\path(0,0)(0,10)(-8.66,15)
    {\thinlines
      \multiput(0,2)(0,1.6){6}{\path(0,0)(-1.73,-1)}
      \multiput(-1.3856,10.8)(-1.3856,0.8){5}{\path(0,0)(-1.73,-1)}}}
   \multiput(86.6,50)(-17.32,0){5}{\path(0,0)(-8.66,-5)(-17.32,0)
    {\thinlines\multiput(-0.866,-0.5)(-1.5588,-0.9){5}{\path(0,0)(0,1.6)}
    \multiput(-8.66,-5)(-1.5588,0.9){6}{\path(0,0)(0,1.6)}}}
% Sector S--+
  {\thinlines
    \multiput(93.8744,45.8)(-1.3856,0.8){6}{\path(0,0)(0,-1.8)}}
  \multiput(86.6,50)(-17.32,0){5}{\path(0,0)(-8.66,5)(-17.32,0)
   {\thinlines\multiput(-0.866,0.5)(-1.5588,0.9){5}{\path(0,0)(0,-1.6)}
    \multiput(-8.66,5)(-1.5588,-0.9){6}{\path(0,0)(0,-1.6)}}}
  \path(95.26,45)(95.26,35)
   {\thinlines
    \multiput(95.26,42)(0,-1.5){6}{\path(0,0)(-1.73,1)}}
  \multiput(95.26,35)(-8.66,-15){3}{\path(0,0)(0,-10)(-8.66,-15)
  {\thinlines\multiput(0,-2)(0,-1.5){6}{\path(0,0)(-1.73,1)}
    \multiput(-0.866,-10.5)(-1.5588,-0.9){6}{\path(0,0)(-1.73,1)}}}
 \put(89,87){$U_1$}
 \put(36,71){$U_2$}
 \put(38,20){$U_3$}
 \put(143,100){$U_1^-$}
 \put(40,100){$U_1^+$}
 \put(98,40){$0$}
\end{picture}
%-----------------------------------
\caption{Sectors of the dual kagome lattice, $d=3$,
$\alpha_k=\exp((2k-1)\pi i/6)$.} \label{fig:sectors}
\end{center}
\end{figure}
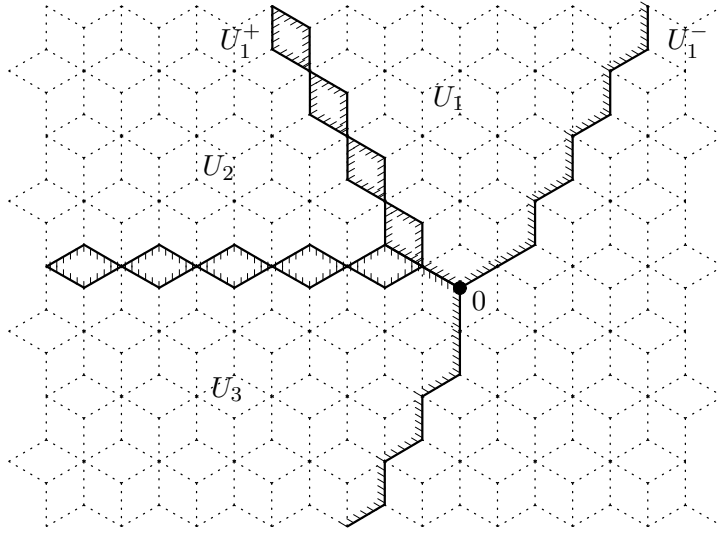

See Fig.\,\ref{fig:sectors} for an illustration.

%-----------------------------------------------------------------

We say that a subset $\cI\subset\Pi_{\ba,\bb}$ is an {\it initial
values locus}, if, prescribing arbitrarily values of $f$ on $\cI$,
one can extend $f$ in virtue of the 3D consistent Cauchy-Riemann
equations in a unique and unambiguous way from $\cI$ to the whole
of $\Pi_{\ba,\bb}$ (cf. \cite{AV}). We will use two types of
initial values loci.
\begin{itemize}
\item
Any monotone path from $\ba$ to $\bb$, with all edges directed
positively:
\begin{equation}\label{IVP2}
\cI_1=\{\bn_r\}_{r=0}^N\quad {\rm with}\quad
\bn_0=\ba,\;\;\bn_N=\bb,\quad {\rm and}\quad
\bn_{r+1}-\bn_r\in\{\be_1,\ldots,\be_d\}.
\end{equation}
\item
The intersection of $\Pi_{\ba,\bb}$ with all coordinate axes:
\begin{equation}\label{IVP1}
\cI_2=\bigcup_{k=1}^d\big\{\bn=n\be_k:\;a_k\le n\le b_k\big\}.
\end{equation}
\end{itemize}
{\bf Proof of Proposition \ref{prop extend}.} We have to show that,
for any set of signs $\beps$, the values of $f$ on 
$V(\Omega_\cD)\cap S_\beps$ determine $f$ uniquely on the hull
$\Pi_{\ba,\bb}\cap S_\beps$.

First, we prove this for $\beps=\beps(m)$, so that $S_\beps=S_m$.
For the sake of notational simplicity, we do this for $m=1$ only,
i.e., for the hull $\cH(V(\Omega_\cD)\cap S_1)=\Pi_{\bO,\bb}$.
Indeed, an arbitrary point $\bn\in V(\Omega_\cD)\cap S_1$ can be
reached from $\bO$ along a path in $V(\Omega_\cD)$ with all edges
from $\{\be_1,\ldots,\be_d\}$. This is a path of the type $\cI_1$,
as in (\ref{IVP2}), hence it is an initial value locus for the
brick $\Pi_{\bO,\bn}$. Since the union of the bricks
$\Pi_{\bO,\bn}$ over all $\bn\in V(\Omega_\cD)\cap S_1$ exhausts
the brick $\Pi_{\bO,\bb}$, our claim is proved.

The bricks $\Pi_{\ba,\bb}\cap S_\beps$ with $\beps\neq \beps(m)$
do not contain points of $V(\Omega_\cD)$ in their interior. However, on
the first step of the proof, we obtain values of $f$ on
all the coordinate axes. This gives an initial values locus of the
type $\cI_2$, as in (\ref{IVP1}), for any brick of this type. \qed
\smallskip

Note that intersections of $\Omega_\cD$ with bricks correspond to
{\em combinatorially convex subsets} of $\cD$, as defined in \cite{M2}.

\section{Discrete exponential functions}
\label{Sect exp}

A particularly important discrete holomorphic function on $\bbZ^d$
is the {\it discrete exponential function}, defined as
\begin{equation}\label{discr exp}
  e(\bn;z)=\prod_{k=1}^d
  \Big(\frac{z+\alpha_k}{z-\alpha_k}\Big)^{n_k}.
\end{equation}
For $d=2$, this function was considered in \cite{F,D1}. The
discrete Cauchy-Riemann equations for the discrete exponential
function are easily checked: they are equivalent to a simple
identity
\[
\Big(\frac{z+\alpha_j}{z-\alpha_j}\cdot
\frac{z+\alpha_k}{z-\alpha_k}-1\Big)\Big/
\Big(\frac{z+\alpha_j}{z-\alpha_j}-
\dfrac{z+\alpha_k}{z-\alpha_k}\Big)=
\frac{\alpha_j+\alpha_k}{\alpha_j-\alpha_k}\,.
\]
At a given $\bn\in\bbZ^d$, the discrete exponential function is
rational with respect to the parameter $z$, with poles at the
points $\epsilon_1\alpha_1,\ldots,\epsilon_d\alpha_d$, where
$\epsilon_k={\rm sign}\, n_k$.

Equivalently, one can identify the discrete exponential function
by its initial values on the axes:
\begin{equation}\label{discr exp ini}
  e(n\be_k;z)=\Big(\frac{z+\alpha_k}{z-\alpha_k}\Big)^n.
\end{equation}
A still another characterization says that $e(\cdot;z)$ is the
B\"acklund transformation of the zero solution of discrete
Cauchy-Riemann equations on $\bbZ^d$, with the ``vertical''
parameter $z$.

Restriction of the function $e(\cdot;z)$ to $V(\Omega_\cD)\sim
V(\cD)$ is a discrete exponential function on $\cD$ defined and
studied in \cite{M1,M2,K}. Note that the latter depends on the
choice of the point $x_0\in V(\cD)$. A question posed in \cite{K}
asks whether discrete exponential functions are dense in the space
of discrete holomorphic functions on $\cD$. We now show that the
answer to this question is in affirmative, in some natural class
of functions (growing not faster than exponentially). 
\begin{theorem}\label{prop dense} 
Let $f$ be a discrete holomorphic function on $V(\cD)\sim V(\Omega_\cD)$, 
satisfying
\begin{equation}\label{exp growth}
  |f(\bn)|\le \exp(C(|n_1|+\ldots+|n_d|)),\qquad \forall \bn\in
  V(\Omega_\cD),
\end{equation}
with some $C\in\bbR$. Extend it to a discrete holomorphic function on
$\cH(V(\Omega_\cD))$. Then inequality (\ref{exp growth}) holds for all 
$\bn\in\cH(V(\Omega_\cD))$, possibly with some larger constant $C$.
There exists a function $g$ defined on the disjoint union of small 
neighborhoods around the points $\pm\alpha_k\in\bbC$ and holomorphic 
on each one of these neighborhoods, such that
\begin{equation}\label{dense}
f(\bn)-f(\bO)=\frac{1}{2\pi i}\int_{\Gamma}
g(\lambda)e(\bn;\lambda)d\lambda, \qquad \forall
\bn\in\cH(V(\Omega_\cD)),
\end{equation}
where $\Gamma$ is a collection of $2d$ small loops, each one
running counterclockwise around one of the points $\;\pm\alpha_k$.
\end{theorem}
{\bf Proof.} In order to extend $f$ from $V(\Omega_\cD)$ to 
$\cH(V(\Omega_\cD))$, one makes elementary steps based on 
eq. (\ref{def multidim holo}). For instance, within the octant
$S_1$ these elementary steps consist of calculating the left-hand side of 
the following equation through the quantities on the right-hand side:
\[
f(\bn+\be_j+\be_k)=f(\bn)+\frac{\alpha_j+\alpha_k}{\alpha_j-\alpha_k}
\big(f(\bn+\be_j)-f(\bn+\be_k)\big).
\]
(In other octants everything is similar, but notations become slightly
more complicated.) A simple induction shows that if the constant $C$ in 
(\ref{exp growth}) satisfies the inequality
\[
1+2\max_{j\neq k}\left|\frac{\alpha_j+\alpha_k}{\alpha_j-\alpha_k}
\right|\exp(C)\le \exp(2C),
\]
then (\ref{exp growth}) propagates in the extension process. This proves
the first statement of the theorem.

To prove the second one, it is enough to find $g(\lambda)$ such
that (\ref{dense}) holds on the coordinate axes, that is,
\begin{equation}\label{dense aux}
f_n^{(k)}-f(\bO)=\underset{\lambda=\alpha_k}{\rm Res}\,
g(\lambda)\left(\frac{\lambda+\alpha_k}{\lambda-\alpha_k}\right)^n,\quad
f_{-n}^{(k)}-f(\bO)=\underset{\lambda=-\alpha_k}{\rm Res}\;
g(\lambda)\left(\frac{\lambda-\alpha_k}{\lambda+\alpha_k}\right)^n,
\quad \forall n>0,
\end{equation}
where $f_n^{(k)}$ are the restrictions of
$f:\cH(V(\Omega_\cD))\to\bbC$ to the coordinate axes:
\[
f_n^{(k)}=f(n\be_k), \qquad a_k(\Omega_\cD)\le n\le
b_k(\Omega_\cD).
\]
Set $g(\lambda)=\sum_{k=1}^d (g_k(\lambda)+g_{-k}(\lambda))$,
where the functions $g_{\pm k}(\lambda)$ vanish everywhere except
in small neighborhoods of the points $\pm\alpha_k$, respectively, and are
given there by convergent series
\begin{equation}\label{dense g_i}
g_k(\lambda)=\frac{1}{2\lambda}\left(f_1^{(k)}-f(\bO)
+\sum_{n=1}^\infty\Big(\frac{\lambda-\alpha_k}{\lambda+\alpha_k}\Big)^n
\big(f_{n+1}^{(k)}-f_{n-1}^{(k)}\big)\right),
\end{equation}
and a similar formula for $g_{-k}(\lambda)$. (Convergence of these 
series is assured by the exponential growth of $f_n^{(k)}$.) The easy-to-check
formula
\[
\underset{\lambda=\alpha_k}{\rm Res}\,\frac{1}{\lambda}
\Big(\frac{\lambda+\alpha_k}{\lambda-\alpha_k}\Big)^n=1-(-1)^n,\quad
n\ge 0,
\]
shows that the so defined function $g$ satisfies (\ref{dense
aux}). \qed

\section{Isomonodromic discrete logarithmic function}
\label{Sect CR isomonodromic}

We first give a construction of the discrete logarithmic function
on $\cD$ which is equivalent to Kenyon's one \cite{K}. This
function is defined, after fixing some point $x_0\in V(\cD)$, by
the formula
\begin{equation}\label{dLog}
 f(x)=\frac{1}{2\pi i}\int_{\Gamma}
 \frac{\log(\lambda)}{2\lambda}\,e(x;\lambda)d\lambda, \qquad
 \forall x\in V(\cD).
\end{equation}
Here the integration path $\Gamma$ is the same as in Theorem
\ref{prop dense}, and fixing $x_0$ is necessary for the definition
of the discrete exponential function on $\cD$. To make
(\ref{dLog}) a valid definition, one has to specify which branch
of $\log(\lambda)$ is chosen around each point $\pm\alpha_k$. This
choice depends on $x$, and is done as follows.

For each $m\in\bbZ$, assign to
$\alpha_m=\exp(i\theta_m)\in\bbS^1$ a certain value of argument
$\theta_m\in\bbR$: choose a value $\theta_1$ of the
argument of $\alpha_1$ arbitrarily, and then extend it according
to the rule
\[
\theta_{m+1}-\theta_m\in (0,\pi),\qquad \forall m\in\bbZ.
\]
Clearly, there holds $\theta_{m+d}=\theta_{m}+\pi$, and therefore
also $\theta_{m+2d}=\theta_{m}+2\pi$. It will be convenient to
consider the points $\alpha_m$, supplied with the arguments
$\theta_m$, as belonging to the Riemann surface
$\widetilde\Lambda$ of the logarithmic function (a branched
covering of the complex $\lambda$-plane). 

The definition domain of the discrete logarithmic function is a
branched covering
\[
 \widetilde U=\bigcup_{m=-\infty}^\infty \widetilde U_m
\]
of the quad-graph $\cD$. Here $\widetilde U_m$ is the sector $U_m$
equipped with additional data -- the interval 
\begin{equation}\label{ineq S}
\log(\alpha_r)\in[i\theta_m,i\theta_{m+d-1}], \qquad r=m,\ldots,m+d-1
\end{equation}
of length less than $\pi$ for the logarithms of the slopes of edges
$\alpha_m,\ldots,\alpha_{m+d-1}$. If $m$ increases by $2d$, the 
interval on the right-hand side of (\ref{ineq S}) is shifted by $2\pi i$.
Two sectors $\widetilde U_{m_1}$
and $\widetilde U_{m_2}$ have a non-empty intersection, if and
only if $|m_1-m_2|<d$. It follows from Lemma \ref{lemma sectors}
that $\widetilde U$ is, indeed, a branched covering of $\cD$.
Definition (\ref{dLog}) should be read as follows: for
$x\in\widetilde U_m$, the poles of $e(x;\lambda)$ are exactly the
points $\alpha_m,\ldots,\alpha_{m+d-1}\in\widetilde\Lambda$. Therefore, 
one can assume that the integration path $\Gamma$ consists of $d$ small 
loops around these points, and the values of $\log(\lambda)$ at these
points satisfy (\ref{ineq S}).
\begin{proposition}{\rm\cite{K}}
The discrete logarithmic function on $\cD$, restricted to
$V(\cG)$, coincides with {\itbf discrete Green's function} on
$\cG$, up to a constant factor $2\pi$.
\end{proposition}
{\bf Proof.} It is not difficult to see that the restriction of
the discrete logarithmic function to black points does not branch:
it is a well-defined real-valued function on  $V(\cG)$. Clearly,
this function is harmonic everywhere except the origin. At the origin,
its Laplacian equals to the increment of $f$ upon running once
around the origin through its white neighbors. The values of $f$
at the vertices neighboring to the origin are nothing but the
arguments of the corresponding edges. Therefore, the above
mentioned increment is equal to $2\pi$. In order to obtain
asymptotic results for the discrete logarithmic function, one can
deform the integration path $\Gamma$ into a connected contour
lying on a single leaf of the Riemann surface of the logarithm,
and then use standard methods of the complex analysis  \cite{K}.
This possibility is due to the fact that functions $g_k$ in
integral representation (\ref{dense}) of an arbitrary discrete
holomorphic function, defined originally in disjoint neighborhoods
of the points $\alpha_r$, in the case of the discrete logarithmic
function are actually restrictions of a single analytic function
$\log(\lambda)/(2\lambda)$ to these neighborhoods. \qed
\medskip

Now we extend the discrete logarithmic function to $\bbZ^d$. To
this end, recall that the sector $U_m$ of $\cD$ is nothing but the
preimage w.r.t. $P$ of the part of $\Omega_\cD$ lying in the
octant $S_m\subset\bbZ^d$. Therefore, it is natural to introduce a
branched covering
\[
\widetilde S=\bigcup_{m=-\infty}^\infty \widetilde S_m
\]
of the set $\bigcup_{m=1}^{2d}S_m\subset\bbZ^d$. Here $\widetilde
S_m$ is the octant $S_m$ equipped with the set of values of
$\log(\epsilon_k\alpha_k)$ satisfying (\ref{ineq S}). Recall that
$\epsilon_k=\epsilon_k(m)$, $k=1,\ldots,d$, are the signs of the
coordinate semi-axes of $S_m$, defined in Sect.\,\ref{Sect CR
multidim}. By definition, $\widetilde S_{m_1}$ and $\widetilde
S_{m_2}$ intersect, if the underlying octants $S_{m_1}$ and
$S_{m_2}$ have a non-empty intersection spanned by the common
coordinate semi-axes, and the data $\log(\epsilon_k\alpha_k)$ for
these common semi-axes match. It is easy to see that $\widetilde
S_{m_1}$ and $\widetilde S_{m_2}$ intersect, if and only if
$|m_1-m_2|<d$.
\begin{definition}
The {\itbf discrete logarithmic function} on $\widetilde S$ is
given by the formula
\begin{equation}\label{Green f prelim}
f(\bn)=\frac{1}{2\pi i}
\int_\Gamma\frac{\log\lambda}{2\lambda}\,e(\bn;\lambda)d\lambda,\qquad
\forall \bn\in\widetilde S,
\end{equation}
where the integration path $\Gamma$ consists, for
$\bn\in\widetilde S_m$, of $d$ loops around
$\alpha_m,\ldots,\alpha_{m+d-1}$, and the branch of the logarithm
on $\Gamma$ is defined by inequality (\ref{ineq S}).
\end{definition}
The discrete logarithmic function on $\cD$ can be described as the
restriction of the discrete logarithmic function on $\widetilde S$
to a branched covering of $\Omega_\cD\sim\cD$. This holds for an
{\em arbitrary} quasicrystallic quad-graph with the set of edge
slopes $A$.

Now we are in a position to give an alternative definition of the
discrete logarithmic function. Clearly, it is completely
characterized by its values $f(n\epsilon_k\be_k)$ on the
coordinate semi-axes of an arbitrary octant $\widetilde S_m$.
\begin{proposition}\label{prop initial values}
For the discrete logarithmic function on $\widetilde S$, each of
$d$ sequences $f_n^{(k)}=f(n\epsilon_k\be_k)$, $k=1,\ldots,d$,
solves the difference equation
\begin{equation}\label{CR axis recur}
 n(f_{n+1}-f_{n-1})=1-(-1)^n,
\end{equation}
with the initial conditions
\begin{equation}\label{init eps}
 f_0^{(k)}=f(\bO)=0,\qquad
 f_1^{(k)}=f(\epsilon_k\be_k)=\log(\epsilon_k\alpha_k).
\end{equation}
Explicitly,
\begin{equation}\label{Green axes}
f_{2n}^{(k)}=\sum_{\ell=1}^n\frac{2}{2\ell-1}\,,\qquad
f_{2n+1}^{(k)}=\log(\epsilon_k\alpha_k),\qquad k=1,\ldots,d,\qquad
n\ge 0.
\end{equation}
Here $\epsilon_k=\epsilon_k(m)$, and the values
$\log(\epsilon_k\alpha_k)$ are chosen in the interval (\ref{ineq
S}).
\end{proposition}
{\bf Proof.} According to eq. (\ref{dense g_i}), the values
$f_n^{(k)}$, with $f_0^{(k)}=0$, are defined by the expansion near
$\lambda=\epsilon_k\alpha_k$,
\begin{equation}\label{Green aux1}
 \log(\lambda)
 =\log(\epsilon_k\alpha_k)+\log\Big(\frac{\lambda}{\epsilon_k\alpha_k}\Big)
 =f_1^{(k)}+\sum_{n=1}^\infty
 \Big(\frac{\lambda-\epsilon_k\alpha_k}{\lambda+\epsilon_k\alpha_k}\Big)^n
 (f_{n+1}^{(k)}-f_{n-1}^{(k)}).
\end{equation}
This is equivalent to
\begin{equation}\label{Green aux2}
 f_1^{(k)}=\log(\epsilon_k\alpha_k),\qquad
 f_{n+1}^{(k)}-f_{n-1}^{(k)}=\frac{1-(-1)^n}{n}.
\end {equation}
The solution to these recurrent relations is given by (\ref{Green
axes}). \qed
\medskip

Observe that values (\ref{Green axes}) at even (resp. odd) points
imitate the behaviour of the real (resp. imaginary) part of the
function $\log(\lambda)$ along the semi-lines ${\rm
arg}(\lambda)={\rm arg}(\epsilon_k\alpha_k)$. This can be easily
extended to the whole of $\widetilde S$. Restricted to black
points $\bn\in\widetilde S$ (those with $n_1+\ldots+n_d$ even),
the discrete logarithmic function models the real part of the
logarithm. In particular, this restricted function is real-valued
and does not branch: its values on $\widetilde S_m$ depend on
$m\pmod{2d}$ only. In other words, it is a well defined function
on $S_m$. On the contrary, the discrete logarithmic function
restricted to white points $\bn\in\widetilde S$ (those with
$n_1+\ldots+n_d$ odd) takes purely imaginary values, and increases
by $2\pi i$, as $m$ increases by $2d$. Hence, this restricted
function models the imaginary part of the logarithm.
\medskip

It turns out that recurrent relations (\ref{CR axis recur})
are characteristic for an important class of solutions of the
discrete Cauchy-Riemann equations, namely for the isomonodromic
ones. Recall the definition of this class. For a discrete
holomorphic function $f:\bbZ^d\to\bbC$, the {\em
transition matrices} are (cf. (\ref{L CR})),
\begin{equation}\label{L CR d}
\renewcommand{\arraystretch}{1.4}
L_k(\bn;\lambda)=\begin{pmatrix} \lambda+\alpha_k &
-2\alpha_k(f(\bn+\be_k)+f(\bn)) \\ 0 & \lambda-\alpha_k
\end{pmatrix}.
\end{equation}
The {\em moving frame} $\Psi(\cdot,\lambda):\bbZ^d\to
GL_2(\bbC)[\lambda]$ is defined by prescribing some $\Psi(\bO;\lambda)$, and
by extending it recurrently according to the formula
\begin{equation}\label{Psi recur}
  \Psi(\bn+\be_k;\lambda)=L_k(\bn;\lambda)\Psi(\bn;\lambda).
\end{equation}
Finally, define the matrices $A(\cdot;\lambda):\bbZ^d\to
GL_2(\bbC)[\lambda]$ by
\begin{equation}\label{A}
  A(\bn;\lambda)=\frac{d\Psi(\bn;\lambda)}{d\lambda}\Psi^{-1}(\bn;\lambda).
\end{equation}
These matrices are defined uniquely after fixing some $A(\bO;\lambda)$.
\begin{definition}
A discrete holomorphic function $f:\bbZ^d\to\bbC$ is called {\itbf
isomonodromic} 
  \footnote{This term originates in the theory of integrable nonlinear 
  differential equations, where it is used for solutions with a similar 
  analytic characterization \cite{IN}.},
if, for some choice of $A(\bO;\lambda)$, the
matrices $A(\bn;\lambda)$ are meromorphic in $\lambda$, with poles
whose positions and orders do not depend on $\bn\in\bbZ^d$.
\end{definition}
It is clear how to extend this definition to functions on the covering
$\widetilde S$.
\begin{theorem}
The discrete logarithmic function is isomonodromic.
\end{theorem}
This is an immediate consequence of the following statement, which we 
formulate for functions on $S_1=(\bbZ_+)^d$ for notational simplicity, 
but which holds actually for any octant $S_\beps$.

\begin{proposition}\label{prop CR isomonodromic}
Let
\begin{equation}\label{A0 CR}
  A(\bO;\lambda)=\frac{1}{\lambda}\begin{pmatrix} 0 & 1 \\
  0 & 0 \end{pmatrix},
\end{equation}
and let there be $d$ sequences $\{f_n^{(k)}\}_{n=0}^\infty$
satisfying, for all $k=1,\ldots,d$, the recurrent relation
(\ref{CR axis recur}). Then the discrete holomorphic function
$f:(\bbZ_+)^d\to\bbC$, defined by the values $f(n\be_k)=f_n^{(k)}$
on the coordinate semi-axes, is isomonodromic. At any point
$\bn\in(\bbZ_+)^d$ there holds:
\begin{equation}\label{A CR simple poles}
  A(\bn;\lambda)=\frac{A^{(0)}(\bn)}{\lambda}+
  \sum_{l=1}^d\Big(\frac{B^{(l)}(\bn)}{\lambda+\alpha_l}+
  \frac{C^{(l)}(\bn)}{\lambda-\alpha_l}\Big),
\end{equation}
with
\begin{eqnarray}
 & A^{(0)}(\bn) =
 \renewcommand{\arraystretch}{1.4}
 \begin{pmatrix}
 0 & (-1)^{n_1+\ldots+n_d} \\ 0 & 0 \end{pmatrix}, &
\label{A CR} \\ \nonumber\\ &
 B^{(l)}(\bn) =
 \renewcommand{\arraystretch}{1.4}
  n_l\begin{pmatrix} 1 & -(f(\bn)+f(\bn-\be_l)) \\
  0 & 0 \end{pmatrix},\qquad
  C^{(l)}(\bn) =
  \renewcommand{\arraystretch}{1.4}
  n_l\begin{pmatrix} 0 & f(\bn+\be_l)+f(\bn) \\
  0 & 1 \end{pmatrix}. & \label{BC CR}
\end{eqnarray}
Moreover, at any point $\bn\in(\bbZ_+)^d$ there holds an
isomonodromic constraint,
\begin{equation}\label{CR constr}
\sum_{l=1}^d n_l\Big(f(\bn+\be_l)-f(\bn-\be_l)\Big)=
1-(-1)^{n_1+\ldots+n_d}.
\end{equation}
\end{proposition}
{\bf Proof.} 
The matrices $A$ satisfy a recurrent relation, which results by
differentiating (\ref{Psi recur}),
\begin{equation}\label{A recur}
  A(\bn+\be_k;\lambda)=
  \frac{dL_k(\bn;\lambda)}{d\lambda}L_k^{-1}(\bn;\lambda)+
  L_k(\bn;\lambda)A(\bn;\lambda)L_k^{-1}(\bn;\lambda).
\end{equation}
Fix some $k=1,\ldots,d$, and consider the
matrices $A(n\be_k;\lambda)$ along the $k$th coordinate semi-axis.
Formula (\ref{A recur}) shows that singularities of
$A(n\be_k;\lambda)$ are poles at $\lambda=0$ and at
$\lambda=\pm\alpha_k$. It is easy to see that the pole $\lambda=0$
remains simple for all $n>0$. As one can show (see Lemma
\ref{lemma CR axes} in Appendix \ref{Sect appendix A}), the
recurrent relation (\ref{CR axis recur}) for $f_n=f(n\be_k)$
assures that the poles $\lambda=\pm\alpha_k$ are simple for all
$n>0$. So, under condition (\ref{CR axis recur}) there holds:
\begin{equation}\label{A CR one pole}
  A(n\be_k;\lambda)=\frac{A^{(0)}(n\be_k)}{\lambda}+
  \frac{B^{(k)}(n\be_k)}{\lambda+\alpha_k}+
  \frac{C^{(k)}(n\be_k)}{\lambda-\alpha_k}\;,
\end{equation}
i.e., eq. (\ref{A CR simple poles}) is valid on the $k$th
coordinate semi-axis, with $B^{(l)}(n\be_k)=C^{(l)}(n\be_k)=0$ for
$l\neq k$. The proof continues by induction, whose scheme follows
filling out the hull of the coordinate semi-axes: each new point
is of the form $\bn+\be_j+\be_k$, $j\neq k$, with three points
$\bn$, $\bn+\be_j$ and $\bn+\be_k$ known from the previous steps,
where the statements of the proposition are assumed to hold. So,
suppose that (\ref{A CR simple poles}) holds at $\bn+\be_j$,
$\bn+\be_k$. The new matrix $A(\bn+\be_j+\be_k;\lambda)$ is
obtained by two alternative formulas,
\begin{equation}\label{prop CR jk}
A(\bn+\be_j+\be_k;\lambda)
 = \bigg(\frac{dL_k(\bn+\be_j;\lambda)}{d\lambda}
 +L_k(\bn+\be_j;\lambda)A(\bn+\be_j;\lambda)\bigg)
 L_k^{-1}(\bn+\be_j;\lambda),
\end{equation}
and the one with interchanged roles of $k$ and $j$. Eq. (\ref{prop
CR jk}) shows that all poles of $A(\bn+\be_j+\be_k;\lambda)$
remain simple, with the possible exception of
$\lambda=\pm\alpha_k$, whose orders might increase by 1. The same
statement holds with $k$ replaced by $j$. Therefore, all poles
remain simple, and (\ref{A CR simple poles}) holds at
$\bn+\be_j+\be_k$. The proof of formulas (\ref{A CR}), (\ref{BC
CR}) and of the constraint (\ref{CR constr}) is based on
computations presented in Appendix \ref{Sect appendix A}. \qed
\medskip

The reason for considering isomonodromic solutions on octants like
$(\bbZ_+)^d$ is clear from (\ref{CR axis recur}): indeed, this
second-order ordinary difference equation has a special form,
enforcing that its solution on the semi-axis $n\ge 0$ is
completely defined by the values at $n=0,1$, and does not depend on 
$f_{-1}$. Thus, a discrete holomorphic function from Proposition 
\ref{prop CR isomonodromic} is uniquely defined by its initial values 
$f(\bO)=f_0$ and $f(\be_k)=f_1^{(k)}$ for $k=1,\ldots,d$.
\medskip

{\bf Remark.} The isomonodromic constraint (\ref{CR constr}) was
found in \cite{NRGO}, without any relation to the discrete logarithmic
function. Observe that our formulation allows us to
avoid a major computational problem arising in \cite{NRGO} in this
context, namely that of compatibility of the constraint with the
discrete Cauchy-Riemann equations.
\medskip

Summing up: {\em discrete Green's function on a quasicrystallic
quad-graph is the real part (i.e., restriction to $V(\cG)$) of the
discrete logarithmic function. The latter can be extended to a
function on a branched covering of certain octants
$S_m\subset\bbZ^d$, $m=1,\ldots,2d$. On each such octant, the
discrete logarithmic function is discrete holomorphic, with the
distinctive property of being isomonodromic. This 
function is uniquely defined either by the integral representation
(\ref{Green f prelim}), or by the values on the coordinate
semi-axes (\ref{Green axes}), or else by the initial values
(\ref{init eps}) and the constraint (\ref{CR constr}).}

\section{3D consistent cross-ratio equations}
\label{Sect cross-ratio}

The cross-ratio system is one of the simplest and at the same time
one of the most fundamental and important {\em nonlinear}
integrable systems on quad-graphs. Recall the definition of the
cross-ratio of four complex numbers:
\begin{equation}\label{q}
q(z_0,z_1,z_2,z_3)=\frac{(z_0-z_1)(z_2-z_3)}{(z_1-z_2)(z_3-z_0)},
\end{equation}
which yields the property
\begin{equation}\label{q sym}
q(z_0,z_1,z_2,z_3)=1/q(z_1,z_2,z_3,z_0).
\end{equation}
Let there be given a function $Q:E(\cG)\sqcup E(\cG^*)\to\bbC$
satisfying the condition
\begin{equation}\label{Q*}
Q(\ee^*)=1/Q(\ee),\qquad\forall \ee\in E(\cG).
\end{equation}
\begin{definition}
A function $z:V(\cD)\to\bbC$ is said to solve the {\itbf
cross--ratio equations} on $\cD$ corresponding to the function
$Q$, if for any positively oriented face $(x_0,y_0,x_1,y_1)$ of
$\cD$ there holds:
\begin{equation}\label{cross-rat eq}
q(z(x_0),z(y_0),z(x_1),z(y_1))=Q(x_0,x_1)=1/Q(y_0,y_1).
\end{equation}
\end{definition}
Like in Sect. \ref{Sect 3D CR}, an interesting question is on the
3D consistency of the system of cross-ratio equations
corresponding to a given function $Q$.

\begin{theorem}\label{Th 3D cr}
The function $Q:E(\cG)\sqcup E(\cG^*)\to\bbC$ can be extended to
$E(\bG)\sqcup E(\bG^*)$ giving a 3D consistent system of
cross-ratio equations, if and only if the following condition is
satisfied:
\begin{equation}\label{cr prod=1}
\prod_{\ee\in{\rm star}(x_0;\cG)} Q(\ee)=1,\quad
 \prod_{\ee^*\in{\rm star}(y_0;\cG^*)} Q(\ee^*)=1,\qquad
 \forall  x_0\in V(\cG),\; y_0\in V(\cG^*).
\end{equation}
\end{theorem}
{\bf Proof.} We proceed as in the proof of Theorem \ref{Th 3D CR}.
Consider a flower of quadrilaterals around $x_0$, with
$\ee_k=(x_0,x_k)$, $\ee_k^*=(y_{k-1},y_k)$. Build the extension of
this flower to the third dimension, as in Sect. \ref{Sect 3D}.
Denote
\begin{equation}\label{Q's}
  Q(x_0,x_k)=Q(\ee_k)=Q_k, \quad Q(x_0,\widehat y_k)=\mu_k.
\end{equation}
Then there holds a statement analogous to Lemma \ref{lemma 3D
nus}: the cross-ratio equations are 3D consistent on the cube over
the $k$-th petal, if and only if
\begin{equation}\label{3D cr aux}
 \mu_{k-1}=Q_k\mu_k.
\end{equation}
This is donstraightforward, as in the proof of Lemma \ref{lemma 3D
nus}. For the cube over the petal with $k=1$, one determines on
the first step the values of $z$ at $x_1$, $\widehat y_0$ and
$\widehat y_1$ from
\begin{eqnarray*}
 q(z(x_0),z(y_0),z(x_1),z(y_1))  & = & Q_1,\\
 q(z(x_0),z(y_0),z(\widehat y_0),z(\widehat x_0)) & = & \mu_0,\\
 q(z(x_0),z(y_1),z(\widehat y_1),z(\widehat x_0)) & = & \mu_1.
\end{eqnarray*}
On the second step one has three alternative values for
$z(\widehat x_1)$ from
\begin{eqnarray*}
 q(z(\widehat x_0),z(\widehat y_0),z(\widehat x_1),z(\widehat y_1))
   & = & Q_1,\\
 q(z(y_1),z(x_1),z(\widehat x_1),z(\widehat y_1)) & = & \mu_0,\\
 q(z(y_0),z(x_1),z(\widehat x_1),z(\widehat y_0)) & = & \mu_1.
\end{eqnarray*}
A direct computation shows that these three values for $z(\widehat
x_1)$ coincide if and only if $\mu_0=Q_1\mu_1$, and then
\[
\widehat x_1=\frac {\mu_0y_0(y_1-\widehat x_0)+\mu_1y_1(\widehat
x_0-y_0)+\widehat x_0(y_0-y_1)}
 {\mu_0(y_1-\widehat x_0)+\mu_1(\widehat x_0-y_0)+(y_0-y_1)}.
\]
Thus, (\ref{3D cr aux}) is proved. This relation yields
immediately that running around $x_0$ returns back an (arbitrary)
initial $\mu_0$, if and only if the first condition in (\ref{cr
prod=1}) holds. The second one follows similarly. \qed

\begin{corollary}
The integrability condition (\ref{cr prod=1}) for the function
$Q:E(\cG)\sqcup E(\cG^*)\to\bbC$ is equivalent to the existence of
a labeling $\alpha^2:E(\cD)\to\bbC$ of undirected edges of $\cD$,
such that, in notations of Fig.\,\ref{diamond again},
\begin{equation}\label{Q paral}
Q(x_0,x_1)=\frac{1}{Q(y_0,y_1)}=\frac{\alpha_0^2}{\alpha_1^2}.
\end{equation}
This formula assures the 3D consistency of the cross-ratio
equations, if one assumes that all vertical edges of $\bD$ carry
one and the same label $\lambda^2\in\bbC$.
\end{corollary}

Let the labeling $\alpha^2$ come from a labeling
$\alpha:\vec{E}(\cD)\to\bbC$ of directed edges. Let
$p:V(\cD)\to\bbC$ be a parallelogram realization of $\cD$ defined
by $p(y)-p(x)=\alpha(x,y)$. Then the cross-ratio equations are
written as
\begin{equation}\label{cr paral}
q(z(x_0),z(y_0),z(x_1),z(y_1))= \frac{\alpha_0^2}{\alpha_1^2}=
q(p(x_0),p(y_0),p(x_1),p(y_1)).
\end{equation}
In other words, the cross-ratio of the vertices of the $f$-image
of any quadrilateral $(x_0,y_0,x_1,y_1)\in F(\cD)$ is equal to the
cross-ratio of the vertices of the corresponding parallelogram.

\begin{proposition}\label{prop cr zcr}
The cross-ratio equations (\ref{cr paral}) admit a zero curvature
representation with the values in $GL_2(\bbC)[\lambda]$, with
transition matrices along $(x,y)\in\vec{E}(\cD)$ given by
\begin{equation}\label{L cr}
\renewcommand{\arraystretch}{1.4}
L(y,x,\alpha;\lambda)=\begin{pmatrix} 1 & z(x)-z(y) \\
\lambda\alpha^2/(z(x)-z(y)) & 1 \end{pmatrix}, \quad where
\quad\alpha=p(y)-p(x).
\end{equation}
\end{proposition}
{\bf Proof.} This result is easy to check. To derive it
systematically using 3D consistency and the procedure outlined at
the end of Sect. \ref{Sect 3D}, observe that eq. (\ref{cr paral})
on the vertical face $(x,y,\widehat y,\widehat x)$ can be written
as
\[
q(z(x),z(y),z(\widehat y),z(\widehat x))= \lambda\alpha^2\quad
\Leftrightarrow \quad z(\widehat y)-z(y)=
L(y,x,\alpha;\lambda)\cdot (z(\widehat x)-z(x)),
\]
where $\lambda=(p(\widehat x)-p(x))^{-2}$, and the matrices
$L(y,x,\alpha;\lambda)$ are as in (\ref{L cr}). One easily shows
that these matrices $L$ form a zero curvature representation with
values in $GL_2(\bbC)[\lambda]$, i.e., that (\ref{zcr prelim})
holds as it stands, and not only up to a scalar factor. \qed
\medskip

The main result of the present section can be formulated as
follows. {\em Integrable cross-ratio equations on a quad-graph
$\cD$ come from parallelogram immersions of $\cD$ in $\bbC$, the
coefficients $Q$ being the cross-ratios of the corresponding
parallelograms. In the case of unitary values $Q\in\bbS^1$ the
parallelograms are actually rhombi.}

\section{Circle patterns and the cross-ratio system}
\label{Sect patterns}

\begin{definition}
A {\itbf Delaunay decomposition} of $\,\bbC$ is a cell decomposition $\cG$
such that the boundary of each face is a polygon inscribed in a circle,
and these circles have no vertices in their interior. These circles build
a {\itbf circle pattern with the combinatorics of} $\cG$.
\end{definition}
%-----------------------------------------------------------------
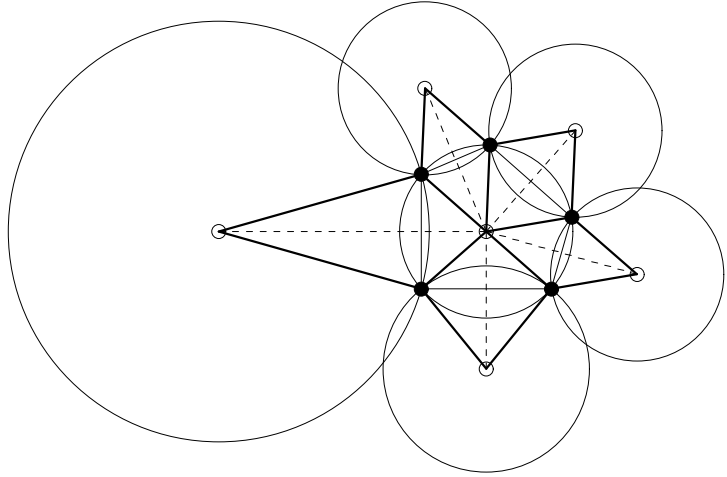
\begin{figure}[htbp]
\begin{center}
\setlength{\unitlength}{0.025em}
\begin{picture}(1000,710)(-700,-360)
 \put(0,0){\circle{20}}
 \put(0,0){\circle{240}}
 \put(5,119.8958){\circle*{20}}
 \put(118.3216,20){\circle*{20}}
 \put(90,-79.3725){\circle*{20}}
 \put(-90,-79.3725){\circle*{20}}
 \put(-90,79.3725){\circle*{20}}
 \put(0,-190){\circle{20}}
 \put(0,-190){\circle{285.2258}}
 \put(208.3216,-59.3725){\circle{20}}
 \put(208.3216,-59.3725){\circle{240}}
 \put(-370,0){\circle{20}}
 \put(-370,0){\circle{582.06527}}
 \put(-85,198.2683){\circle{20}}
 \put(-85,198.2683){\circle{240}}
 \put(123.3216,139.8958){\circle{20}}
 \put(123.3216,139.8958){\circle{240}}
 \path(5,119.8958)(118.3216,20)
 \path(118.3216,20)(90,-79.3725)
 \path(90,-79.3725)(-90,-79.3725)
 \path(-90,-79.3725)(-90,79.3725)
 \path(-90,79.3725)(5,119.8958)
 \dashline[+30]{10}(0,0)(0,-190)
 \dashline[+30]{10}(0,0)(208.3216,-59.3725)
 \dashline[+30]{10}(0,0)(-370,0)
 \dashline[+30]{10}(0,0)(-85,198.2683)
 \dashline[+30]{10}(0,0)(123.3216,139.8958)
 \thicklines
 \path(-370,0)(-90,-79.3725)  \path(-370,0)(-90,79.3725)
 \path(0,0)(-90,-79.3725)  \path(0,0)(-90,79.3725)
 \path(-85,198.2683)(-90,79.3725)  \path(-85,198.2683)(5,119.8958)
 \path(0,0)(5,119.8958)
 \path(123.3216,139.8958)(5,119.8958)  \path(123.3216,139.8958)(118.3216,20)
 \path(0,0)(118.3216,20)
 \path(208.3216,-59.3725)(118.3216,20) \path(208.3216,-59.3725)(90,-79.3725)
 \path(0,0)(90,-79.3725)
 \path(0,-190)(90,-79.3725)  \path(0,-190)(-90,-79.3725)
    \end{picture}
\caption{Circle pattern}\label{circle pattern}
\end{center}
\end{figure}
%-----------------------------------------------------------------
The vertices $z:V(\cG)\to\bbC$ of a Delaunay decomposition are
the intersection points of the circles of the corresponding pattern.
The circle of the pattern corresponding to a face $y\in F(\cG)$ will
be denoted by $C(y)$. If two faces $y_0,y_1\in F(\cG)$ 
have a common edge $(x_0,x_1)$, then the circles $C(y_0)$
and $C(y_1)$ intersect in the points $z(x_0),z(x_1)$. In other
words, the edges of $\cG$ correspond to pairs of neighboring
(intersecting) circles of the pattern. Similarly, if several faces
$y_1,y_2,\ldots,y_m$ of $\cG$ meet at one vertex $x_0\in V(\cG)$,
then the corresponding circles $C(y_1),C(y_2),\ldots,C(y_m)$ also
have a common intersection point $z(x_0)$.

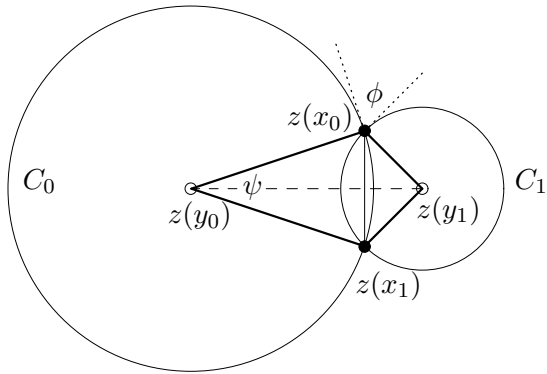
\begin{figure}[htbp]
\begin{center}
\setlength{\unitlength}{0.04em}
\begin{picture}(500,340)(-300,-160)
 \put(-150,0){\circle{10}}
 \put(-170,-30){$z(y_0)$}
 \put(50,0){\circle{10}}
 \put(45,-25){$z(y_1)$}
 \put(0,50){\circle*{10}}
 \put(-67,55){$z(x_0)$}
 \put(0,-50){\circle*{10}}
 \put(-8,-90){$z(x_1)$}
 \put(50,0){\circle{141.42136}}
 \put(130,0){$C_1$}
 \put(-150,0){\circle{316.22777}}
 \put(-295,0){$C_0$}
 \path(0,50)(0,-50)
 \dashline[+30]{9}(-150,0)(-110,0)
 \dashline[+30]{10}(-85,0)(50,0)
 \dottedline{6}(0,50)(50,100)
 \dottedline{6}(0,50)(-25,125)
 \thicklines
 \path(-150,0)(0,50)  \path(0,50)(50,0)
 \path(50,0)(0,-50)  \path(0,-50)(-150,0)
  \put(1,75){$\phi$}
  \put(-105,-5){$\psi$}
    \end{picture}
\caption{Two intersecting circles}\label{two circles}
\end{center}
\end{figure}

Given a circle pattern with the combinatorics of $\cG$, we can
extend the function $z$ to the vertices of the dual graph $\cG^*$,
setting
\[
z(y)={\rm center\;\;of\;\;the\;\;circle}\;\;C(y), \quad y\in
F(\cG)\simeq V(\cG^*).
\]
After this extension, the map $z$ is defined on all of
$V(\cD)=V(\cG)\sqcup V(\cG^*)$, where $\cD$ is the double of $\cG$.
Consider a face of the double. Its vertices $x_0,x_1,y_0,y_1$
correspond to the intersection points and to the centers of two
neighboring circles $C_0,C_1$ of the pattern. The following statement
is obtained by a simple computation.
\begin{lemma}\label{lemma on angles}
If $\phi$ is the intersection angle of $C_0,C_1$, as on
Fig.\,\ref{two circles}, then
\begin{equation}
q(z(x_0),z(y_0),z(x_1),z(y_1))=\exp(2i\phi).
\end{equation}
\end{lemma}
%-----------------------------------------------------------------
%-----------------------------------------------------------------
It will be convenient to assign the intersection angle $\phi$ of
$C(y_0),C(y_1)$ to the edge $(y_0,y_1)\in E(\cG^*)$. Extend the
function $\phi:E(\cG^*)\to(0,\pi)$ to $E(\cG)$ by setting
$\phi(\ee)=\pi-\phi(\ee^*)$.

\begin{proposition}\label{labelling for pattern}
Let $\cG$ be Delaunay decomposition of a plane, and
consider a circle pattern with the combinatorics of $\cG$ and with
the intersection angles $\phi:E(\cG^*)\to(0,\pi)$. Let
$\{z(x):x\in V(\cG)\}$ and $\{z(y):y\in V(\cG^*)\}$ consist of 
intersection points of the circles, resp. of their
centers. Then $z:V(\cD)\to\Bbb C$ satisfies a system of
cross-ratio equations with the function $Q:E(\cG)\sqcup
E(\cG^*)\to\bbS^1$ defined as $Q(\ee)=\exp(2i\phi(\ee))$. There
holds:
\begin{equation}\label{CP angles int points}
  \prod_{\ee\in{\rm star}(x_0;\cG)}\exp(2i\phi(\ee))=1,
  \quad \forall x_0\in V(\cG).
\end{equation}
The following condition is necessary and sufficient for the
integrability of the system of cross-ratio equations:
\begin{equation}\label{CP angles centers}
  \prod_{\ee^*\in{\rm star}(y_0;\cG^*)}\exp(2i\phi(\ee^*))=1,
  \quad \forall y_0\in V(\cG^*),
\end{equation}
i.e., for each circle of the pattern the sum of its intersection
angles with all neighboring circles of the pattern vanishes
$\pmod\pi$.
\end{proposition}
{\bf Proof.} The relation (\ref{CP angles int points}) is obvious
for geometrical reasons: for an arbitrary common intersection
point of circles of the pattern, the sum of their consecutive
pairwise intersection angles vanishes $\pmod\pi$. Now the claim
follows from Theorem \ref{Th 3D cr}. \qed
\medskip

We can formulate the main result of this section as follows. {\em
Combinatorial data $\cG$ and intersection angles
$\phi:E(\cG^*)\to(0,\pi)$ belong to an integrable circle pattern,
if and only if they admit an isoradial realization. This latter
realization gives a rhombic immersion of the double $\cD$, and
generates also a dual isoradial circle pattern with the
combinatorial data $\cG^*$ and intersection angles
$\phi:E(\cG)\to(0,\pi)$.}

\section{Hirota system}
\label{Sect Hirota}

We have seen that integrable circle patterns deliver solutions
$z:V(\cD)\to\bbC$ to integrable cross-ratio systems. These
solutions are characterized as follows: the $z$-image of any
quadrilateral $(x_0,y_0,x_1,y_1)$ from $F(\cD)$ is a {\it kite}
with two pairs of sides of equal length (incident to the white
vertices $z(y_0)$, $z(y_1)$), and with the prescribed angle
$\pi-\phi$ at the black vertices $z(x_0)$, $z(x_1)$. The following 
transformation of the cross-ratio system is useful in order to single 
out this class of kite solutions.
\begin{definition}
Let $\alpha:\vec{E}(\cD)\to\bbC$ be a labeling, and let
$p:V(\cD)\to\bbC$ be the corresponding parallelogram realization
of $\cD$ defined by $p(y)-p(x)=\alpha(x,y)$. A function
$w:V(\cD)\to\bbC$ is said to solve the corresponding {\itbf Hirota
system}, if for any positively oriented face $(x_0,y_0,x_1,y_1)\in
F(\cD)$ there holds, in the notations of Fig.\,\ref{diamond again
or}:
\begin{equation}\label{Hirota}
\alpha_0w(x_0)w(y_0)+\alpha_1w(y_0)w(x_1)-\alpha_0w(x_1)w(y_1)-
\alpha_1w(y_1)w(x_0)=0,
\end{equation}
or, in a more invariant fashion,
\begin{eqnarray}\label{Hirota in p}
\lefteqn{w(x_0)w(y_0)\big(p(y_0)-p(x_0)\big)+
         w(y_0)w(x_1)\big(p(x_1)-p(y_0)\big)+}\nonumber\\
      && w(x_1)w(y_1)\big(p(y_1)-p(x_1)\big)+
         w(y_1)w(x_0)\big(p(x_0)-p(y_1)\big)=0.
\end{eqnarray}
\end{definition}
Obviously, a black-white scaling, i.e., a transformation $w\to cw$
on $V(\cG)$ and $w\to c^{-1}w$ on $V(\cG^*)$ with a constant $c$,
maps solutions of the Hirota system into solutions. We will
identify solutions related by such a transformation.

\begin{proposition}
Let $w:V(\cD)\to\bbC$ be a solution of the Hirota system. Then the
relation
\begin{equation}\label{w vs z}
z(y)-z(x)=\alpha(x,y) w(x)w(y)=w(x)w(y)\big(p(y)-p(x)\big), \qquad
\forall (x,y)\in\vec{E}(\cD),
\end{equation}
correctly defines a unique (up to an additive constant) function
$z:V(\cD)\to\bbC$ which is a solution of the cross-ratio system
(\ref{cr paral}). Conversely, for any solution $z$ of the
cross-ratio system (\ref{cr paral}), relation (\ref{w vs z})
defines a function $w$ correctly and uniquely (up to a black-white
scaling); this function $w$ solves the Hirota system
(\ref{Hirota}).
\end{proposition}
{\bf Proof.} Simple calculation based on closing conditions around
the quadrilateral $(x_0,y_0,x_1,y_1)$. \qed

\begin{proposition}
Let all $\alpha\in\bbS^1$, so that $p:V(\cD)\to\bbC$ is a rhombic
realization of $\cD$. Let $z:V(\cD)\to\bbC$ be a solution of the
corresponding cross-ratio system (\ref{cr paral}). It corresponds
to a circle pattern, if and only if the corresponding function $w$
satisfies the condition
\begin{equation}\label{Hirota red}
w(x)\in\bbS^1,\quad w(y)\in\bbR_+,\quad \forall x\in V(\cG),\;y\in
V(\cG^*).
\end{equation}
The values $w:V(\cG^*)\to\bbR_+$ have then the interpretation of
the radii of the circles.
\end{proposition}
{\bf Proof.} The function $z$ corresponds to a circle pattern, if
and only if all elementary quadrilaterals
$(z(x_0),z(y_0),z(x_1),z(y_1))$ are of the kite form with the
properties:

-- the pairs of edges incident with white vertices have equal
length,

-- the angles at black vertices are equal to the corresponding
angles of the underlying rhombi.

\noindent As easily seen, these conditions are equivalent to:
\[
\frac{|w(x_0)|}{|w(x_1)|}=1 \quad{\rm and}\quad
\frac{w(y_0)}{w(y_1)}\in\bbR_+,
\]
respectively. This yields (\ref{Hirota red}), possibly upon a
black-white scaling. \qed
\smallskip

{\bf Remark.} The conditions (\ref{Hirota red}) form an {\em
admissible reduction} of the Hirota system corresponding to a
rhombic realization, in the following sense: if any three of the
four points $w(x_0)$, $w(y_0)$, $w(x_1)$, $w(y_1)$ satisfy the
condition (\ref{Hirota red}), then so does the fourth one. This is
immediately seen, if one rewrites the Hirota equation
(\ref{Hirota}) in one of the equivalent forms:
\begin{equation}\label{Hirota alt}
\frac{w(x_1)}{w(x_0)}=
\frac{\alpha_1w(y_1)-\alpha_0w(y_0)}{\alpha_1w(y_0)-\alpha_0w(y_1)}\quad
\Leftrightarrow\quad \frac{w(y_1)}{w(y_0)}=
\frac{\alpha_0w(x_0)+\alpha_1w(x_1)}{\alpha_0w(x_1)+\alpha_1w(x_0)}.
\end{equation}

\begin{proposition}
a) Let $\alpha:\vec{E}(\cD)\to\bbC$ be a labeling, and let
$p:V(\cD)\to\bbC$ be the corresponding parallelogram realization
of $\cD$ defined by $p(y)-p(x)=\alpha(x,y)$. Then the
corresponding Hirota system is 3D consistent.

b) Let all $\alpha\in\bbS^1$, so that $p:V(\cD)\to\bbC$ is a
rhombic realization of $\cD$. Consider a solution
$w:V(\cD)\to\bbC$ corresponding to a circle pattern with the
combinatorics of $\cG$, i.e., satisfying (\ref{Hirota red}).
Consider its B\"acklund transformation $\widehat w:V(\cD)\to\bbC$
with an arbitrary parameter $\lambda\in\bbS^1$ and with an
arbitrary initial value $\widehat w(x_0)\in\bbR_+$ or $\widehat
w(y_0)\in\bbS^1$. Then there holds:
\begin{equation}\label{Hirota red hat}
\widehat w(x)\in\bbR_+,\quad \widehat w(y)\in\bbS^1,\quad \forall
x\in V(\cG),\;y\in V(\cG^*),
\end{equation}
so that $\widehat w$ corresponds to a circle pattern with the
combinatorics of $\cG^*$.
\end{proposition}
{\bf Proof.} Statement a) is a matter of a direct computation. In
the notations of Fig.\,\ref{cube}, suppose that the Hirota
equation (\ref{Hirota}) holds on all faces of the cube, wherein
the vertical edges carry the (arbitrary) label $\lambda$.  One
finds that all three alternative ways to compute $w(\widehat x_1)$
lead to one and the same result, namely
\[
w(\widehat x_1)=\frac{\lambda(\alpha_0^2-\alpha_1^2)w(y_0)w(y_1)+
\alpha_1(\lambda^2-\alpha_0^2)w(y_0)w(\widehat x_0)+
\alpha_0(\alpha_1^2-\lambda^2)w(y_1)w(\widehat x_0)}
{\lambda(\alpha_0^2-\alpha_1^2)w(\widehat x_0)+
\alpha_1(\lambda^2-\alpha_0^2)w(y_1)+
\alpha_0(\alpha_1^2-\lambda^2)w(y_0)}.
\]
Statement b) follows from the Remark above. \qed

\begin{proposition}\label{prop Hirota zcr}
The Hirota system (\ref{Hirota}) admits a zero curvature
representation with the values in $GL_2(\bbC)[\lambda]$, with
transition matrices along $(x,y)\in\vec{E}(\cD)$ given by
\begin{equation}\label{L Hirota}
\renewcommand{\arraystretch}{1.4}
L(y,x,\alpha;\lambda)=\begin{pmatrix} 1 & -\alpha w(y) \\
-\lambda\alpha/w(x) & w(y)/w(x) \end{pmatrix}, \quad where
\quad\alpha=p(y)-p(x).
\end{equation}
\end{proposition}
{\bf Proof.}  The matrix (\ref{L Hirota}) is obtained also
directly from (\ref{L cr}) by the substitution (\ref{w vs z}),
followed by a simple gauge transformation
\[
L\mapsto \begin{pmatrix} 1 & 0 \\ 0 & w(y) \end{pmatrix}L
\begin{pmatrix} 1 & 0 \\ 0 & 1/w(x) \end{pmatrix}.
\]
Alternatively, a systematic derivation of this result is based on
the 3D consistency and the procedure outlined at the end of Sect.
\ref{Sect 3D}. \qed
\medskip

The main result of the present section is as follows. {\em
Integrable circle patterns can be alternatively described by
solutions of the Hirota system with a special property of being
real-valued on $V(\cG)$ and unimodular on $V(\cG^*)$.}

\section{Linearization}
\label{Sect linearization}

Let $\alpha:\vec{E}(\cD)\to\bbC$ be a labeling, and let
$p:V(\cD)\to\bbC$ be the corresponding parallelogram realization
of $\cD$ defined by $p(y)-p(x)=\alpha(x,y)$. Then the formula
\[
z_0(x)=p(x),\qquad w_0(x)=1,\qquad\forall x\in V(\cD),
\]
gives a (trivial) solution of the cross-ratio system (\ref{cr
paral}) and the corresponding (trivial) solution of the Hirota
system. Suppose that $z_0:V(\cD)\to\bbC$ belongs to a differentiable
one-parameter family of solutions $z_\epsilon:V(\cD)\to\bbC$, $\epsilon\in
(-\epsilon_0,\epsilon_0)$, of the same cross-ratio system, and
denote by $w_\epsilon:V(\cD)\to\bbC$
the corresponding solutions of the Hirota system. Denote
\begin{equation}\label{linearization}
g=\frac{dz_{\epsilon}}{d\epsilon}\bigg|_{\epsilon=0}\;,\qquad
f=\bigg(w_\epsilon^{-1}\,\frac{dw_{\epsilon}}{d\epsilon}\bigg)_{\epsilon=0}\;.
\end{equation}
\begin{proposition}\label{prop linearization}
Both functions $f,g:V(\cD)\to\bbC$ solve discrete Cauchy-Riemann equations
(\ref{CR paral}).
\end{proposition}
{\bf Proof.} By differentiating (\ref{w vs z}), we obtain a relation 
between the functions $f,g:V(\cD)\to\bbC$:
\begin{equation}\label{f vs g}
g(y)-g(x)=\big(f(x)+f(y)\big)\big(p(y)-p(x)\big), \qquad \forall
(x,y)\in\vec{E}(\cD).
\end{equation}
The proof of proposition is based on this relation solely. Indeed,
the closeness condition for the form on the right-hand side reads:
\begin{eqnarray*}
\lefteqn{\big(f(x_0)+f(y_0)\big)\big(p(y_0)-p(x_0)\big)+
         \big(f(y_0)+f(x_1)\big)\big(p(x_1)-p(y_0)\big)+}\nonumber\\
      && \big(f(x_1)+f(y_1)\big)\big(p(y_1)-p(x_1)\big)+
         \big(f(y_1)+f(x_0)\big)\big(p(x_0)-p(y_1)\big)=0,
\end{eqnarray*}
which is equivalent to (\ref{CR paral}) for the function $f$.
Similarly, the closeness condition for $f$, that is,
\[
\big(f(x_0)+f(y_0)\big)-\big(f(y_0)+f(x_1)\big)+
\big(f(x_1)+f(y_1)\big)-\big(f(y_1)+f(x_0)\big)=0,
\]
yields:
\[
\frac{g(y_0)-g(x_0)}{p(y_0)-p(x_0)}-\frac{g(x_1)-g(y_0)}{p(x_1)-p(y_0)}
+\frac{g(y_1)-g(x_1)}{p(y_1)-p(x_1)}-\frac{g(x_0)-g(y_1)}{p(x_0)-p(y_1)}=0.
\]
Under the condition $p(y_0)-p(x_0)=p(x_1)-p(y_1)$, this is
equivalent to (\ref{CR paral}) for $g$.   \qed
\medskip

{\bf Remark.} 
This proof shows that, given a discrete holomorphic
function $f:V(\cD)\to\bbC$, relation (\ref{f vs g}) correctly 
defines a unique, up to an additive constant, function
$g:V(\cD)\to\bbC$, which is also discrete holomorphic. Conversely, 
for any $g$ satisfying the discrete Cauchy-Riemann equations 
(\ref{CR paral}), relation (\ref{f vs g}) defines a function $f$ 
correctly and uniquely (up to an additive black-white constant); this 
function $f$ also solves the discrete Cauchy-Riemann equations 
(\ref{CR paral}). Actually, formula (\ref{f vs g}) expresses 
that the discrete holomorphic function $f$ is the {\it discrete derivative}
of $g$, so that $g$ is obtained from $f$ by {\it discrete
integration}. This operation was considered in \cite{D1, D2, M1}.
\medskip

Summarizing, we have the following statement.
\begin{theorem}
a) A tangent space to the set of solutions of an integrable
cross-ratio system, at a point corresponding to a rhombic embedding of 
a quad-graph, consists of discrete holomorphic functions on this
embedding. This holds in both descriptions of the above set: in terms 
of variables $z$ satisfying the cross-ratio equations, and in
terms of variables $w$ satisfying the Hirota equations. The 
corresponding two descriptions of the tangent space are related by 
taking the discrete derivative (resp. anti-derivative) of discrete 
holomorphic functions.  

b) A tangent space to the set of integrable circle
patterns of a given combinatorics, at a point corresponding to
an isoradial pattern, consists of discrete holomorphic functions
on the rhombic embedding of the corresponding quad-graph, which
take real values at white vertices and purely imaginary values at black 
ones. This holds in the description of circle patterns in terms of 
circle radii and rotation angles at intersection points (Hirota equations).
\end{theorem}

\section{Isomonodromic discrete power function}
\label{Sect z^c}

Like in Sect. \ref{Sect CR multidim}, one can consider functions
$z:\bbZ^d\to\bbC$ and $w:\bbZ^d\to\bbC$, satisfying, on each
elementary square of $\bbZ^d$, the cross-ratio and the Hirota
equation, respectively, and ask about isomonodromic solutions. As
shown in \cite{AB1, AB2, BH}, this leads to a discrete analog of
the power function. Since the latter references contain a detailed
presentation of these results in terms of the cross-ratio
variables $z$, we restrict ourselves here to similar results in
terms of the Hirota variables $w$. (Recall that transition
matrices in these two formulations actually coincide, up to a
simple gauge transformation by diagonal matrices which do not
depend on $\lambda$.)

Transition matrices for the Hirota system on $\bbZ^d$ are:
\begin{equation}\label{L Hirota d}
\renewcommand{\arraystretch}{1.4}
L_k(\bn;\lambda)= \begin{pmatrix} 1 & -\alpha_k w(\bn+\be_k) \\
-\lambda\alpha_k/w(\bn) & w(\bn+\be_k)/w(\bn) \end{pmatrix}.
\end{equation}
With these transition matrices, isomonodromic solutions are
defined in exactly the same manner as in Sect. \ref{Sect CR
isomonodromic}.

\begin{proposition}\label{prop Hirota isomonodromic}
Let
\begin{equation}\label{A0 Hirota}
  A(\bO;\lambda)=\frac{1}{\lambda}\begin{pmatrix} -\gamma/2 & 0 \\
  0 & \gamma/2 \end{pmatrix},
\end{equation}
and let there be $d$ sequences $\{w_n^{(k)}\}_{n=0}^\infty$
satisfying, for all $k=1,\ldots,d$, the recurrent relation
\begin{equation}\label{Hirota axis recur}
 n\,\frac{w_{n+1}-w_{n-1}}{w_{n+1}+w_{n-1}}=\big(\gamma-\tfrac{1}{2}\big)
 \big(1-(-1)^n\big).
\end{equation}
Then the solution $w:(\bbZ_+)^d\to\bbC$ of the Hirota system,
defined by the values $w(n\be_k)=w_n^{(k)}$ on the coordinate
semi-axes, is isomonodromic. At any point $\bn\in(\bbZ_+)^d$ there
holds:
\begin{equation}\label{A Hirota simple poles}
  A(\bn;\lambda)=\frac{A^{(0)}(\bn)}{\lambda}+
  \sum_{l=1}^d\frac{B^{(l)}(\bn)}{\lambda-\alpha_l^{-2}}\;,
\end{equation}
with
\begin{eqnarray}
 & A^{(0)}(\bn) =
 \renewcommand{\arraystretch}{1.4}
 \begin{pmatrix}
 -\gamma/2 & * \\ 0 & \gamma/2 \end{pmatrix}, &
\label{A Hirota} \\ \nonumber\\ &
 B^{(l)}(\bn) =
 \renewcommand{\arraystretch}{1.4}
  \dfrac{n_l}{w(\bn+\be_l)+w(\bn-\be_l)}
  \begin{pmatrix} w(\bn+\be_l) & \alpha_l w(\bn+\be_l)w(\bn-\be_l)\\
  1/\alpha_l & w(\bn-\be_l) \end{pmatrix}. & \label{B Hirota}
\end{eqnarray}
The upper right entry of the matrix $A^{(0)}(\bn)$, denoted by the
asterisk in (\ref{A Hirota}), is actually given by
\begin{equation}\label{A12 Hirota}
A_{12}^{(0)}(\bn)=-\sum_{l=1}^d B_{12}^{(l)}(\bn).
\end{equation}
Moreover, at any point $\bn\in(\bbZ_+)^d$ there holds an
isomonodromic constraint,
\begin{equation}\label{Hirota constr}
\sum_{l=1}^dn_l\,\frac{w(\bn+\be_l)-w(\bn-\be_l)}
 {w(\bn+\be_l)+w(\bn-\be_l)}=
 \big(\gamma-\tfrac{1}{2}\big)\big(1-(-1)^{n_1+\ldots+n_d}\Big).
\end{equation}
\end{proposition}
{\bf Proof.} The scheme of the proof is the same as for
Proposition \ref{prop CR isomonodromic}. Fix some $k=1,\ldots,d$,
and consider the matrices $A(n\be_k;\lambda)$ along the $k$th
coordinate semi-axis. It follows from formula (\ref{A recur}) that
singularities of $A(n\be_k;\lambda)$ are poles at $\lambda=0$ and
at $\lambda=\alpha_k^{-2}$. While the pole $\lambda=0$
automatically remains simple for all $n>0$, this is not
necessarily so for the pole $\lambda=\alpha_k^{-2}$. As one can
show (see Lemma \ref{lemma Hirota axes} in Appendix \ref{Sect
appendix B}), the recurrent relation (\ref{Hirota axis recur}) for
$w_n=f(n\be_k)$ assures that the pole $\lambda=\alpha_k^{-2}$ is
simple for all $n>0$, and
\begin{equation}\label{A Hirota one pole}
  A(n\be_k;\lambda)=\frac{A^{(0)}(n\be_k)}{\lambda}+
  \frac{B^{(k)}(n\be_k)}{\lambda-\alpha_k^{-2}}\;,
\end{equation}
i.e., eq. (\ref{A Hirota simple poles}) is valid on the $k$th
coordinate semi-axis, with $B^{(l)}(n\be_k)=0$ for $l\neq k$. To
prove that eq. (\ref{A Hirota simple poles}) is valid also
elsewhere, one argues by induction: suppose that (\ref{A Hirota
simple poles}) holds at $\bn+\be_j$, $\bn+\be_k$. Then eq.
(\ref{prop CR jk}) shows that all poles of
$A(\bn+\be_j+\be_k;\lambda)$ remain simple, with the possible
exception of $\lambda=\alpha_k^{-2}$, whose order might increase
by 1. The same statement holds with $k$ replaced by $j$. Hence,
all poles remain simple. Therefore, (\ref{A Hirota simple poles})
holds at $\bn+\be_j+\be_k$, possibly up to a term which does not
vanish with $\lambda\to\infty$. Such a term is absent, if the
right-hand side of (\ref{A recur}) vanishes with
$\lambda\to\infty$, that is, if
\begin{equation}\label{A12 Hirota alt}
\begin{pmatrix} 0 & 0 \\ 1 & 0 \end{pmatrix}
\bigg(A^{(0)}(\bn)+\sum_{l=1}^dB^{(l)}(\bn)\bigg)
\begin{pmatrix} 0 & 0 \\ 1 & 0 \end{pmatrix}=0.
\end{equation}
Clearly, the latter equation is equivalent to (\ref{A12 Hirota}).
Computations towards the proof of (\ref{A12 Hirota}), as well as
of (\ref{A Hirota}), (\ref{B Hirota}) and of constraint
(\ref{Hirota constr}), are presented in Appendix \ref{Sect
appendix B}. \qed
\medskip

{\bf Remark.} Again, the isomonodromic constraint (\ref{Hirota
constr}) was found in \cite{NRGO}. In the approach of that paper,
consistency of the constraint with the Hirota equations (called
lattice MKdV there) is a difficult problem, only manageable with
the help of a computer system for symbolic manipulations. In our
formulation, this comes for free, as a natural consequence of the
construction based on the 3D consistency of the Hirota system.
\medskip

A solution of the Hirota system given in Proposition \ref{prop
Hirota isomonodromic} is completely defined by its initial values
$w(0)=w_0$ and $w(\be_k)=w_1^{(k)}$ for $k=1,\ldots,d$. The choice
\begin{equation}\label{Hirota init}
  w_0=1,\qquad w_1^{(k)}=\exp(i\rho_k)\,, \quad k=1,\ldots,d,
\end{equation}
with arbitrary constants $\rho_k$, leads to the following solution
on the semi-axes:
\begin{equation}\label{Hirota axes}
w_{2n}^{(k)}=\prod_{\ell=1}^n\frac{\ell-1+\gamma}{\ell-\gamma}\,,\qquad
w_{2n+1}^{(k)}=\exp(i\rho_k).
\end{equation}
Observe the asymptotics at $n\to\infty$,
\begin{equation}
w_{2n}^{(k)}=c(\gamma)n^{2\gamma-1}\big(1+O(n^{-1})\big).
\end{equation}
The following special choice of $\rho_k$ defines the {\it discrete
analog of the function $w\mapsto w^{2\gamma-1}$ on} $(\bbZ_+)^d$:
\begin{equation}\label{power w1}
i\rho_k=(2\gamma-1)\log\alpha_k, \quad{\rm so\quad that}\quad
w(\be_k)=\alpha_k^{2\gamma-1}.
\end{equation}

{\bf Remark.} In the variables $z$ the initial values
$\{z_n^{(k)}\}_{n=0}^\infty$ on the semi-axes are given by the
following analog (and consequence) of (\ref{Hirota axis recur}):
\begin{equation}\label{cr axis recur}
 n\,\frac{(z_{n+1}-z_n)(z_n-z_{n-1})}{z_{n+1}-z_{n-1}}=\gamma z_n,
\end{equation}
and then the corresponding solution of the cross-ratio equation
satisfies an analog of (\ref{Hirota constr}):
\begin{equation}\label{cr constr}
\sum_{j=1}^d
n_j\,\frac{(z(\bn+\be_j)-z(\bn))(z(\bn)-z(\bn-\be_j))}
{z(\bn+\be_j)-z(\bn-\be_j)}=\gamma z(\bn).
\end{equation}
Recall also \cite{AB1, AB2, BH} that the discrete analog of the
function $z\mapsto z^{2\gamma}$ on $(\bbZ_+)^d$ is characterized
by the constraint (\ref{cr axis recur}) and the following choice
of the initial conditions:
\begin{equation}\label{cr init}
z(0)=0,\qquad z(\be_k)=\alpha_k^{2\gamma},\quad k=1,\ldots,d.
\end{equation}
Clearly, this choice of initial conditions is equivalent to
(\ref{Hirota init}), if one takes into account the basic relation
(\ref{w vs z}) between the variables $z$ and $w$.
\medskip

Like in Sect. \ref{Sect CR isomonodromic}, isomonodromic solutions
similar to those of Proposition \ref{prop Hirota isomonodromic}
can be defined not only on $(\bbZ_+)^d$ but on any octant
$S_\beps$. They are characterized by the initial data
\begin{equation}\label{Hirota init eps}
w(0)=0,\qquad
w(\epsilon_k\be_k)=(\epsilon_k\alpha_k)^{2\gamma-1},\quad
k=1,\ldots,d,
\end{equation}
and give discrete analogs of the function $w\mapsto w^{2\gamma-1}$
on $S_\beps$. Such a solution is fixed by an independent choice of
branches of the function $w^{2\gamma-1}$ at the points
$w=\epsilon_k\alpha_k$. This is equivalent to choosing the
branches of the function $\log w$, because of
$w^{2\gamma-1}=\exp((2\gamma-1)\log w)$.
\begin{definition}
The {\itbf discrete power function} $w^{2\gamma-1}$ on $\widetilde
S$ is a complex-valued function whose restriction to $\widetilde
S_m$ is defined as the unique isomonodromic solution
$w:S_m\to\bbC$ of the Hirota system on the corresponding $S_m$
with the initial data (\ref{Hirota init eps}) fixed by the
condition (\ref{ineq S}).
\end{definition}
Clearly, the discrete power function takes real values at the
white points and unimodular values at the black points, so that it
corresponds to a circle pattern.
\begin{proposition}
The tangent vector to the space of integrable circle patterns
along the curve consisting of patterns $w^{2\gamma-1}$, at the
point corresponding to $\gamma=1/2$, is the discrete logarithmic
function.
\end{proposition}
{\bf Proof.} We have to prove that the discrete logarithm $f$ and
the discrete power function $w^{2\gamma-1}$ are related by
\[
f(\bn)=\Big(\frac{1}{2}\frac{d}{d\gamma}
\,w^{2\gamma-1}(\bn)\Big)_{\gamma=1/2}\,.
\]
It is enough to prove this on the coordinate semi-axes of each
octant $S_m$. But this follows by differentiating with respect to
$\gamma$ constraint (\ref{Hirota axis recur}) and initial
conditions (\ref{Hirota init eps}) at the point $\gamma=1/2$,
where all $w=1$: the results coincide with (\ref{CR axis recur})
and (\ref{init eps}), respectively. \qed

\section{Concluding remarks}
\label{Sect conclusions}

Results of Sect. \ref{Sect linearization} can be generalized to
the case of linearization at an arbitrary (not necessarily
parallelogram) solution $z_0$ of the cross-ratio system and the
corresponding solution $w_0$ of the Hirota system. In this case
the relation between derivatives (\ref{linearization}), coming to
replace eq. (\ref{f vs g}), reads:
\begin{equation}\label{f vs g nonrhombic}
g(y)-g(x)=\big(f(x)+f(y)\big)\big(z_0(y)-z_0(x)\big), \qquad
\forall (x,y)\in\vec{E}(\cD).
\end{equation}
Arguments similar to those of the proof of Proposition \ref{prop
linearization} show that in this case the function $f$ is discrete
holomorphic with respect to $z_0$, i.e.,
\[
\frac{f(y_1)-f(y_0)}{f(x_1)-f(x_0)}=
\frac{z_0(y_1)-z_0(y_0)}{z_0(x_1)-z_0(x_0)}\,, \qquad \forall
(x_0,y_0,x_1,y_1)\in F(\cD).
\]
The function $g$, in general, is no longer discrete holomorphic.

Thus, a tangent space to the set of integrable circle patterns of
a given combinatorics, at an arbitrary point, consists of
functions, discrete holomorphic with respect to the kite-form
embedding $z_0$ of the corresponding quad-graph. This holds for
the description of circle patterns in terms of circle radii and
rotation angles at intersection points (Hirota equations). The
elements of the tangent space are characterized by the property of
being real at white vertices and purely imaginary at black ones.

A number of constructions of the present paper can be generalized
to the case of kite-form (rather than rhombic) embeddings coming
from an integrable circle pattern. In particular, differentiating
the discrete $w^{2\gamma-1}$ with respect to $\gamma$ at a point
$\gamma\neq 1/2$, one obtains a sort of the discrete logarithmic
(and Green's) functions on the kite-form quad-graph corresponding
to $z^{2\gamma}$.

\begin{appendix}
\section{Appendix: proof of Proposition \ref{prop CR isomonodromic}}
\label{Sect appendix A}

\begin{lemma}\label{lemma CR axes}
Let the matrix $A(\bO;\lambda)$ be as in (\ref{A0 CR}). Fix some
$k=1,\ldots,d$. Then singularities of the matrices
$A(n\be_k;\lambda)$ are poles at $\lambda=0$,
$\lambda=\pm\alpha_k$. For $n>0$, the poles $\lambda=0$ and
$\lambda=-\alpha_k$ are simple. The pole $\lambda=\alpha_k$ is
simple for all $n>0$, if and only if recurrent relation (\ref{CR
axis recur}) holds for $f_n=f(n\be_k)$. In this case there holds
(\ref{A CR one pole}) with
\begin{eqnarray}
 & A^{(0)}(n\be_k) =
 \renewcommand{\arraystretch}{1.4}
 \begin{pmatrix} 0 & (-1)^n \\ 0 & 0 \end{pmatrix}, &
\label{A CR axis}\\ \nonumber\\
 & B^{(k)}(n\be_k) =
  \renewcommand{\arraystretch}{1.4}
  n\begin{pmatrix} 1 & -(f_n+f_{n-1}) \\
  0 & 0 \end{pmatrix},\qquad
  C^{(k)}(n\be_k) =
  \renewcommand{\arraystretch}{1.4}
  n\begin{pmatrix} 0 & f_{n+1}+f_n \\
  0 & 1 \end{pmatrix}. & \label{BC CR axis}
\end{eqnarray}
\end{lemma}
{\bf Proof} proceeds by induction. Putting matrices (\ref{L CR d})
into the recurrent definition (\ref{A recur}), one finds
immediately:
\[
A_{11}(n\be_k;\lambda)=\frac{n}{\lambda+\alpha_k}\,,\qquad
A_{22}(n\be_k;\lambda)=\frac{n}{\lambda-\alpha_k}\,,
\]
and the following recurrent relation for the upper right entry of
the matrix $A(\cdot;\lambda)$,
\[
A_{12}((n+1)\be_k;\lambda) =
\frac{\lambda+\alpha_k}{\lambda-\alpha_k}\, A_{12}(n\be_k;\lambda)
+\frac{2\alpha_k(f_{n+1}+f_n)}{\lambda-\alpha_k}\left(
\frac{n+1}{\lambda+\alpha_k}-\frac{n}{\lambda-\alpha_k}\right).
\]
Assume that there holds eq. (\ref{A CR one pole}). Then it holds
also with $n\mapsto n+1$, if and only if no higher order pole
appears at $\lambda=\alpha_k$ by this transition, what is
equivalent to
\begin{equation}\label{A CR lemma constr}
C_{12}^{(k)}(n\be_k)=n(f_{n+1}+f_n).
\end{equation}
This has to be considered as a recursive definition of the
sequence $\{f_n\}$ (isomonodromic constraint). Under this
condition, we find the recurrent relations for the upper right
entries of the matrices $A^{(0)}$, $B^{(k)}$, $C^{(k)}$:
\begin{eqnarray}
A_{12}^{(0)}((n+1)\be_k) & = & -A_{12}^{(0)}(n\be_k), \label{A CR
lemma a recur}\\ B_{12}^{(k)}((n+1)\be_k) & = &
-(n+1)(f_{n+1}+f_n), \label{A CR lemma b recur}\\
C_{12}^{(k)}((n+1)\be_k) & = & (n+1)(f_{n+1}+f_n)+
B_{12}^{(k)}(n\be_k)+C_{12}^{(k)}(n\be_k)+2A_{12}^{(0)}(n\be_k).
\label{A CR lemma c recur}
\end{eqnarray}
Now (\ref{A CR lemma a recur}), (\ref{A CR lemma b recur}) yield:
\begin{equation}\label{A CR lemma ab}
A_{12}^{(0)}(n\be_k)=(-1)^n, \qquad
B_{12}^{(k)}(n\be_k)=-n(f_n+f_{n-1}).
\end{equation}
Adding all three equations (\ref{A CR lemma a recur})--(\ref{A CR
lemma c recur}), we find that
\begin{equation}
A_{12}^{(0)}(n\be_k)+B_{12}^{(k)}(n\be_k)+C_{12}^{(k)}(n\be_k)=1,
\end{equation}
and this together with (\ref{A CR lemma constr}), (\ref{A CR lemma
ab}) implies explicit form (\ref{CR axis recur}) of the
isomonodromic constraint. \qed
\medskip

{\bf Proof of Proposition \ref{prop CR isomonodromic}, continued.}
As shown in the main text, the induction from $\bn+\be_j$,
$\bn+\be_k$ to $\bn+\be_j+\be_k$ proves formula (\ref{A CR simple
poles}). From the diagonal part of eq. (\ref{prop CR jk}) one
easily derives that for all $\bn\in\bbZ^d$,
\[
A_{11}(\bn;\lambda)=\sum_{l=1}^d\frac{n_l}{\lambda+\alpha_l}\,,\qquad
A_{22}(\bn;\lambda)=\sum_{l=1}^d\frac{n_l}{\lambda-\alpha_l}\,.
\]
The following formula is an easy consequence of eq. (\ref{prop CR
jk}) under the limit $\lambda\to\infty$:
\begin{equation}\label{CR prop constr}
A_{12}^{(0)}(\bn)+\sum_{l=1}^d
\Big(B_{12}^{(l)}(\bn)+C_{12}^{(l)}(\bn)\Big)=1.
\end{equation}
It remains to show that the following relations propagate in the
evolution defined by the recurrent relation (\ref{prop CR jk}):
\begin{eqnarray}
 & A_{12}^{(0)}(\bn)=(-1)^{n_1+\ldots+n_d}\,, & \label{CR prop A12}\\
 & B_{12}^{(l)}(\bn)=-n_l(f(\bn)+f(\bn-\be_l)),\qquad
C_{12}^{(l)}(\bn)=n_l(f(\bn+\be_l)+f(\bn)). & \label{CR prop BC12}
\end{eqnarray}
Indeed, constraint (\ref{CR constr}) follows then immediately,
because it coincides with (\ref{CR prop constr}), if one takes
(\ref{CR prop A12})--(\ref{CR prop BC12}) into account. Writing
now the upper right entry of eq. (\ref{prop CR jk}) in length, we
find the following recurrent relations:
\begin{eqnarray}
A_{12}^{(0)}(\bn+\be_j+\be_k) & = & -A_{12}^{(0)}(\bn+\be_j),
 \label{CR prop a recur}\\  \nonumber\\
B_{12}^{(l)}(\bn+\be_j+\be_k) & \!\! = & \!\!
 \frac{\alpha_l-\alpha_k}{\alpha_l+\alpha_k}\,B_{12}^{(l)}(\bn+\be_j)-
 \frac{2\alpha_k}{\alpha_l+\alpha_k}(n_l+\delta_{lj}+\delta_{lk})
 (f(\bn+\be_j+\be_k)+f(\bn+\be_j)), \nonumber \\
 \label{CR prop b recur} \\
C_{12}^{(l)}(\bn+\be_j+\be_k) & = &
 \frac{\alpha_l+\alpha_k}{\alpha_l-\alpha_k}\,C_{12}^{(l)}(\bn+\be_j)-
 \frac{2\alpha_k}{\alpha_l-\alpha_k}(n_l+\delta_{lj})
 (f(\bn+\be_j+\be_k)+f(\bn+\be_j)),\nonumber\\
 \label{CR prop ck recur}
\end{eqnarray}
the latter formula being valid for $l\neq k$ only. For
$C_{12}^{(k)}(\bn+\be_j+\be_k)$ there holds a similar but much
longer formula, which we actually will not need. Now, eq. (\ref{CR
prop a recur}) readily yields (\ref{CR prop A12}). By the
induction hypothesis, eqs. (\ref{CR prop b recur}) and (\ref{CR
prop ck recur}) with $l\neq k$ can be rewritten as
\begin{eqnarray*}
\lefteqn{B_{12}^{(l)}(\bn+\be_j+\be_k)=}\\
 &&-(n_l+\delta_{lj})\bigg(\frac{\alpha_l-\alpha_k}{\alpha_l+\alpha_k}
 (f(\bn+\be_j)+f(\bn+\be_j-\be_l))
 +\frac{2\alpha_k}{\alpha_l+\alpha_k}(f(\bn+\be_j+\be_k)+f(\bn+\be_j))
 \bigg),\\
\lefteqn{C_{12}^{(l)}(\bn+\be_j+\be_k)=}\\
 &&(n_l+\delta_{lj})\bigg(\frac{\alpha_l+\alpha_k}{\alpha_l-\alpha_k}
 (f(\bn+\be_j+\be_l)+f(\bn+\be_j))
 -\frac{2\alpha_k}{\alpha_l-\alpha_k}(f(\bn+\be_j+\be_k)+f(\bn+\be_j))
 \bigg).
\end{eqnarray*}
But the discrete Cauchy-Riemann equation for the corresponding
elementary squares imply that the latter two equations are
equivalent to
\begin{eqnarray*}
B_{12}^{(l)}(\bn+\be_j+\be_k) & = & -(n_l+\delta_{lj})
(f(\bn+\be_j+\be_k)+f(\bn+\be_j+\be_k-\be_l)),\\
C_{12}^{(l)}(\bn+\be_j+\be_k) & = & (n_l+\delta_{lj})
(f(\bn+\be_j+\be_k+\be_l)+f(\bn+\be_j+\be_k)),
\end{eqnarray*}
which coincide with (\ref{CR prop BC12}) at $\bn+\be_j+\be_k$ for
$l\neq k$. By interchanging the roles of $k$ and $j$, we see that
(\ref{CR prop BC12}) at $\bn+\be_j+\be_k$ holds also for $l\neq
j$, and thus for all $l=1,\ldots,d$. \qed

\section{Appendix: proof of Proposition \ref{prop Hirota isomonodromic}}
\label{Sect appendix B}

\begin{lemma}\label{lemma Hirota axes}
Let the matrix $A(\bO;\lambda)$ be as in (\ref{A0 Hirota}). Fix some
$k=1,\ldots,d$. Then singularities of the matrices
$A(n\be_k;\lambda)$ are poles at $\lambda=0$,
$\lambda=\alpha_k^{-2}$. The pole $\lambda=0$ is simple. The pole
$\lambda=\alpha_k^{-2}$ is simple for all $n>0$, if recurrent
relation (\ref{cr axis recur}) holds for $w_n=w(n\be_k)$. In this
case there holds (\ref{A Hirota one pole}) with
\begin{eqnarray}
 & A^{(0)}(n\be_k) =
 \renewcommand{\arraystretch}{1.4}
 \begin{pmatrix} -\gamma/2 & * \\ 0 & \gamma/2 \end{pmatrix}, &
\label{A Hirota axis}\\ \nonumber\\
 & B^{(k)}(n\be_k) =
  \renewcommand{\arraystretch}{1.4}
  \dfrac{n}{w_{n+1}+w_{n-1}}
  \begin{pmatrix} w_{n+1} & \alpha_k w_{n+1}w_{n-1} \\
  1/\alpha_k & w_{n-1} \end{pmatrix}. & \label{B Hirota axis}
\end{eqnarray}
The upper right entry of the matrix $A^{(0)}(n\be_k)$, denoted in
(\ref{A Hirota axis}) by the asterisk, is given by
\begin{equation}\label{A12 Hirota axis}
A_{12}^{(0)}(n\be_k)=-B_{12}^{(k)}(n\be_k).
\end{equation}
\end{lemma}
{\bf Proof.} Assume that eq. (\ref{A Hirota one pole}) holds. Put
matrices (\ref{L Hirota d}) into recurrent definition (\ref{A
recur}). It is easy to see that no higher order pole appears at
$\lambda=\alpha_k^{-2}$ by the transition $n\mapsto n+1$, if and
only if
\[
\renewcommand{\arraystretch}{1.4}
\begin{pmatrix} 1 & -\alpha_kw_{n+1} \\
-1/(\alpha_kw_n) & w_{n+1}/w_n \end{pmatrix} B^{(k)}(n\be_k)
\begin{pmatrix} 1 & \alpha_kw_n  \\
1/(\alpha_kw_{n+1}) & w_n/w_{n+1} \end{pmatrix}=0.
\]
This is equivalent to
\begin{equation}\label{lemma Hirota isom B}
\renewcommand{\arraystretch}{1.4}
\begin{pmatrix} 1 & -\alpha_kw_{n+1}\end{pmatrix} B^{(k)}(n\be_k)
\begin{pmatrix} 1 \\  1/(\alpha_kw_{n+1}) \end{pmatrix}=0,
\end{equation}
or, written in length, to
\begin{equation}\label{lemma Hirota isom B again}
B_{22}^{(k)}(n\be_k)-B_{11}^{(k)}(n\be_k) +\alpha_kw_{n+1}
B_{21}^{(k)}(n\be_k)
-\frac{1}{\alpha_kw_{n+1}}B_{12}^{(k)}(n\be_k)=0.
\end{equation}
This is a recursive definition of the sequence $\{w_n\}$ (an
isomonodromic constraint). Notice that this is a {\em quadratic}
equation for $w_{n+1}$, unlike (\ref{A CR lemma constr}), which
was a {\em linear} equation for $f_{n+1}$.

Under condition (\ref{lemma Hirota isom B}), or (\ref{lemma Hirota
isom B again}), eq. (\ref{A Hirota one pole}) holds also by
$n\mapsto n+1$, possibly with an additional $\lambda$-independent
term on the right-hand side, which vanishes if and only if
(\ref{A12 Hirota axis}) holds. We will show in a moment that this
is indeed the case. One readily finds recurrent relations for the
matrices $A^{(0)}(n\be_k)$ and $B^{(k)}(n\be_k)$. For the matrix
$A^{(0)}(n\be_k)$ they read:
\begin{equation}\label{lemma Hirota isom A}
\renewcommand{\arraystretch}{1.4}
A^{(0)}((n+1)\be_k)= \begin{pmatrix} 1 & -\alpha_kw_{n+1} \\ 0 &
w_{n+1}/w_n \end{pmatrix} A^{(0)}(n\be_k)
\begin{pmatrix} 1 & \alpha_kw_n  \\ 0 & w_n/w_{n+1} \end{pmatrix}.
\end{equation}
This proves formula (\ref{A Hirota axis}), with a recurrent
relation for the upper right entry:
\begin{equation}\label{lemma Hirota isom A12}
A_{12}^{(0)}((n+1)\be_k)= -\gamma\alpha_kw_n+
\frac{w_n}{w_{n+1}}A_{12}^{(0)}(n\be_k).
\end{equation}
For the matrix $B^{(k)}(n\be_k)$ the recurrent relations read, in
components:
\begin{eqnarray}
B_{11}^{(k)}((n+1)\be_k) & = &
  \gamma+B_{22}^{(k)}(n\be_k)-
  \frac{1}{\alpha_kw_{n+1}}\big(A_{12}^{(0)}(n\be_k)+
  B_{12}^{(k)}(n\be_k)\big),  \label{lemma Hirota isom B11}\\
B_{22}^{(k)}((n+1)\be_k) & = &
  1-\gamma+B_{11}^{(k)}(n\be_k)+\frac{1}{\alpha_kw_{n+1}}
  \big(A_{12}^{(0)}(n\be_k)+B_{12}^{(k)}(n\be_k)\big),
  \label{lemma Hirota isom B22}\\
B_{12}^{(k)}((n+1)\be_k) & = & \gamma\alpha_kw_n-
  \frac{w_n}{w_{n+1}}A_{12}^{(0)}(n\be_k),
  \label{lemma Hirota isom B12}\\
B_{21}^{(k)}((n+1)\be_k) & = & \frac{1-\gamma}{\alpha_kw_n}+
  \frac{1}{\alpha_k^2w_nw_{n+1}}
  \big(A_{12}^{(0)}(n\be_k)+2B_{12}^{(k)}(n\be_k)\big)\nonumber\\
&& -\frac{1}{\alpha_kw_n}
  \big(B_{22}^{(k)}(n\be_k)-B_{11}^{(k)}(n\be_k)\big).
  \label{lemma Hirota isom B21}
\end{eqnarray}
Comparing now (\ref{lemma Hirota isom B12}) with (\ref{lemma
Hirota isom A12}), we see that (\ref{A12 Hirota axis}) holds for
all $n$, as claimed above. Upon using this fact and constraint
(\ref{lemma Hirota isom B again}), we can rewrite formulas
(\ref{lemma Hirota isom B11})--(\ref{lemma Hirota isom B21}) as
follows:
\begin{eqnarray}
B_{11}^{(k)}((n+1)\be_k) & = &
  \gamma+B_{22}^{(k)}(n\be_k),
  \label{lemma Hirota isom B11 again}\\
B_{22}^{(k)}((n+1)\be_k) & = &
  1-\gamma+B_{11}^{(k)}(n\be_k),
  \label{lemma Hirota isom B22 again}\\
B_{12}^{(k)}((n+1)\be_k) & = & \gamma\alpha_kw_n+
  \frac{w_n}{w_{n+1}}B_{12}^{(k)}(n\be_k),
  \label{lemma Hirota isom B12 again}\\
B_{21}^{(k)}((n+1)\be_k) & = & \frac{1-\gamma}{\alpha_kw_n}+
  \frac{w_{n+1}}{w_n} \,B_{21}^{(k)}(n\be_k).
  \label{lemma Hirota isom B21 again}
\end{eqnarray}
These relations together with constraint (\ref{lemma Hirota isom B
again}) define the sequence $\{w_n\}$ and the matrices
$B^{(k)}(n\be_k)$ completely. First of all, there follows from
(\ref{lemma Hirota isom B11 again}), (\ref{lemma Hirota isom B22
again}):
\begin{equation}\label{lemma Hirota isom B1122}
B_{11}^{(k)}(n\be_k)+B_{22}^{(k)}(n\be_k)=n,\qquad
B_{11}^{(k)}(n\be_k)-B_{22}^{(k)}(n\be_k)=\big(\gamma-\tfrac{1}{2}\big)
\big(1-(-1)^n\big).
\end{equation}
Further, there follows from (\ref{lemma Hirota isom B11
again})--(\ref{lemma Hirota isom B21 again}):
\begin{eqnarray}
B_{11}^{(k)}((n+1)\be_k)-\frac{1}{\alpha_kw_n}B_{12}^{(k)}((n+1)\be_k)
& = &
  B_{22}^{(k)}(n\be_k)-\frac{1}{\alpha_kw_{n+1}}B_{12}^{(k)}(n\be_k),
  \label{lemma Hirota isom B1112}\\
B_{22}^{(k)}((n+1)\be_k)-\alpha_kw_nB_{21}^{(k)}((n+1)\be_k) & = &
  B_{11}^{(k)}(n\be_k)-\alpha_kw_{n+1}B_{21}^{(k)}(n\be_k).
  \label{lemma Hirota isom B2221}
\end{eqnarray}
Subtracting these two equations and taking (\ref{lemma Hirota isom
B again}) into account, we find, upon downshifting $n$:
\begin{equation}\label{lemma Hirota isom B shifted}
B_{11}^{(k)}(n\be_k)-B_{22}^{(k)}(n\be_k)
+\alpha_kw_{n-1}B_{21}^{(k)}(n\be_k)
-\frac{1}{\alpha_kw_{n-1}}B_{12}^{(k)}(n\be_k)=0.
\end{equation}
This yields that one of the solutions of eq. (\ref{lemma Hirota
isom B again}), considered as a quadratic equation for $w_{n+1}$,
is $w_{n+1}=-w_{n-1}$. We will be interested in the second one. To
find it, add eqs. (\ref{lemma Hirota isom B again}) and
(\ref{lemma Hirota isom B shifted}) and derive, under the
condition $w_{n+1}+w_{n-1}\neq 0$:
\begin{equation}\label{lemma Hirota isom B1221}
B_{12}^{(k)}(n\be_k)=\alpha_k^2w_{n+1}w_{n-1}B_{21}^{(k)}(n\be_k).
\end{equation}
Due to (\ref{lemma Hirota isom B again}), the right--hand sides of
(\ref{lemma Hirota isom B1112}), (\ref{lemma Hirota isom B2221})
are equal to one another. Using there (\ref{lemma Hirota isom
B1221}), we finally come to
\begin{equation}\label{lemma Hirota isom B11 and 22}
B_{11}^{(k)}(n\be_k)=\alpha_kw_{n+1}B_{21}^{(k)}(n\be_k),\qquad
B_{22}^{(k)}(n\be_k)=\alpha_kw_{n-1}B_{21}^{(k)}(n\be_k).
\end{equation}
This together with (\ref{lemma Hirota isom B1122}) yields
\begin{equation}\label{lemma Hirota isom B21 answer}
B_{21}^{(k)}(n\be_k)=\frac{n}{\alpha_k(w_{n+1}+w_{n-1})},
\end{equation}
and now both the expression (\ref{B Hirota axis}) and explicit
form (\ref{Hirota axis recur}) of the constraint (\ref{lemma
Hirota isom B again}) follow readily. \qed
\medskip

{\bf Proof of Proposition \ref{prop Hirota isomonodromic},
continued.} As demonstrated in the main text, if eq. (\ref{A
Hirota simple poles}) holds at $\bn+\be_j$, $\bn+\be_k$, then it
holds also at $\bn+\be_j+\be_k$, provided eq. (\ref{A12 Hirota})
is valid. To prove eq. (\ref{A12 Hirota}), put expression (\ref{L
Hirota d}) into eq. (\ref{prop CR jk}), and find recurrent
relations for the matrices $A^{(0)}(\bn)$ and $B^{(l)}(\bn)$. Upon
use of the abbreviation $L_{j,k}(\lambda)=L_k(\bn+\be_j;\lambda)$,
we have:
\begin{eqnarray}
A^{(0)}(\bn+\be_j+\be_k) & = & L_{j,k}(0)A^{(0)}(\bn+\be_j)
 L_{j+k,-k}(0), \label{prop Hirota A recur}\\
B^{(l)}(\bn+\be_j+\be_k) & = &
 L_{j,k}(\alpha_l^{-2})\,
 \dfrac{B^{(l)}(\bn+\be_j)}{1-\alpha_k^2\alpha_l^{-2}}\,
 L_{j+k,-k}(\alpha_l^{-2}), \quad l\neq k,\qquad
 \label{prop Hirota B recur}\\
B^{(k)}(\bn+\be_j+\be_k) & = & -L_{j,k}(\alpha_k^{-2})
 \bigg(A^{(0)}(n+\be_j)+\sum_{l\neq k}
 \frac{B^{(l)}(\bn+\be_j)}{1-\alpha_k^2\alpha_l^{-2}}\bigg)
 L_{j+k,-k}(\alpha_k^{-2})\nonumber\\
& & + \;\;{\rm lower\;\;triangular\;\;matrix}.
 \label{prop Hirota Bk recur}
\end{eqnarray}
Here we used the fact that
$L_{j,k}^{-1}(\lambda)=L_{j+k,-k}(\lambda)/(1-\lambda\alpha_k^{-2})$.
Taking into account that the upper triangular part of the matrix
$L_{j,k}(\lambda)$ does not depend on $\lambda$ (and coincides
with $L_{j,k}(0)$), we see that eqs. (\ref{prop Hirota A
recur})--(\ref{prop Hirota Bk recur}) imply the desired property
(\ref{A12 Hirota}), which proves formula (\ref{A Hirota simple
poles}) at $\bn+\be_j+\be_k$.

After eq. (\ref{A Hirota simple poles}) has been proved, it is
instructive to rewrite eq. (\ref{A recur}) as
\[
  A(\bn+\be_k;\lambda)L_k(\bn;\lambda)-L_k(\bn;\lambda)A(\bn;\lambda)=
  \frac{dL_k(\bn;\lambda)}{d\lambda},
\]
and consider the limit $\lambda\to\infty$ of {\em this} formula.
Due to (\ref{L Hirota d}) and (\ref{A Hirota simple poles}), this
limit reads:
\begin{equation}\label{prop Hirota matrix constraint}
\bigg(A^{(0)}(\bn+\be_k)+\sum_{l=1}^dB^{(l)}(\bn+\be_k)\bigg)
\begin{pmatrix} 0 & 0 \\ 1 & 0 \end{pmatrix}-
\begin{pmatrix} 0 & 0 \\ 1 & 0 \end{pmatrix}
\bigg(A^{(0)}(\bn)+\sum_{l=1}^dB^{(l)}(\bn)\bigg)=
\begin{pmatrix} 0 & 0 \\ 1 & 0 \end{pmatrix}.
\end{equation}
Clearly, this equation contains more information than (\ref{A12
Hirota alt}). More precisely, the diagonal terms of (\ref{prop
Hirota matrix constraint}) are equivalent to (\ref{A12 Hirota
alt}), while the lower left entry gives an additional identity.

Further, eq. (\ref{prop Hirota A recur}) yields immediately that
the matrix $A^{(0)}(\bn+\be_j+\be_k)$ retains the upper triangular
form (\ref{A Hirota}). Turning to the matrices $B^{(l)}(\bn)$,
observe first of all that formula (\ref{B Hirota}) is equivalent
to
\begin{equation}\label{B Hirota alt}
 B^{(l)}(\bn) \sim
 \renewcommand{\arraystretch}{1.4}
  \begin{pmatrix} \alpha_lw(\bn+\be_l) \\ 1 \end{pmatrix}
  \begin{pmatrix} 1 & \alpha_l w(\bn-\be_l) \end{pmatrix},\qquad
  {\rm tr}\, B^{(l)}(\bn)=n_l,
\end{equation}
where the sign $\sim$ means ``equal up to a scalar factor''.
Suppose that this holds at the points $\bn+\be_j$ and $\bn+\be_k$
for all $l=1,\ldots,d$. Then it follows from (\ref{prop Hirota B
recur}) that for all $l\neq k$ the matrix
$B^{(l)}(\bn+\be_j+\be_k)$ is also of rank 1, and its trace is
equal to ${\rm tr}\,
B^{(l)}(\bn+\be_j)=n_l+\delta_{jl}=n_l+\delta_{jl}+\delta_{kl}$.
It remains to prove that
\begin{eqnarray}
& L_{j,k}(\alpha_l^{-2})\begin{pmatrix} \alpha_lw(\bn+\be_j+\be_l)
 \\ 1 \end{pmatrix} \sim \begin{pmatrix} \alpha_lw(\bn+\be_j+\be_k+\be_l)
 \\ 1 \end{pmatrix}, & \\ \nonumber\\
& \begin{pmatrix} 1 & \alpha_l w(\bn+\be_j-\be_l)\end{pmatrix}
 L_{j+k,-k}(\alpha_l^{-2}) \sim
  \begin{pmatrix} 1 & \alpha_l
  w(\bn+\be_j+\be_k-\be_l)\end{pmatrix}. &
\end{eqnarray}
These equations, written in length, read:
\begin{eqnarray}
 w(\bn+\be_j)\,\dfrac{\alpha_kw(\bn+\be_j+\be_k)-\alpha_lw(\bn+\be_j+\be_l)}
 {\alpha_kw(\bn+\be_j+\be_l)-\alpha_lw(\bn+\be_j+\be_k)}
 & = & w(\bn+\be_j+\be_k+\be_l), \\
 w(\bn+\be_j)\,\dfrac{\alpha_kw(\bn+\be_j+\be_k)+\alpha_lw(\bn+\be_j-\be_l)}
 {\alpha_kw(\bn+\be_j-\be_l)+\alpha_lw(\bn+\be_j+\be_k)}
 & = & w(\bn+\be_j+\be_k-\be_l),
\end{eqnarray}
and are nothing but the Hirota equations on the corresponding
elementary squares. Thus, recurrent relation (\ref{prop Hirota B
recur}) implies that formula (\ref{B Hirota alt}) holds at
$\bn+\be_j+\be_k$ for $l\neq k$. By interchanging the roles of $j$
and $k$, formula (\ref{B Hirota alt}) holds for $l\neq j$, and
thus for all $l=1,\ldots,d$.
It remains to prove the isomonodromic constraint (\ref{Hirota
constr}). But it is not difficult to see that it is a direct
consequence of the lower left entry of eq. (\ref{prop Hirota
matrix constraint}), if one takes into account expressions (\ref{A
Hirota}) and (\ref{B Hirota}). \qed
\end{appendix}
%%%%%%%%%%%%%%%%%%%%%%%%%%%%%%%%%%%%%%%%%%%%%%%%%%%%%%%%%%%%%%%%%%%%%%%%%%

\end{document}